\documentclass[12pt,smallextended]{svjour3}
\smartqed
\usepackage{lineno}
\modulolinenumbers[5]

\journalname{Journal of XXX}
\usepackage{amssymb,amsmath,amsfonts}
\usepackage[bookmarksnumbered, colorlinks, urlcolor=blue, linkcolor=blue, citecolor=red, plainpages=false, pdfstartview=FitW]{hyperref}
\usepackage{array}
\usepackage{bm}
\usepackage{epsfig}
\usepackage{mathptmx}
\usepackage{diagbox}
\usepackage{graphicx}
\usepackage{color}
\usepackage{caption}
\usepackage{wrapfig}
\usepackage{graphics}
\usepackage{enumerate}
\usepackage{amsxtra}
\usepackage{latexsym}
\usepackage{multirow}
\usepackage{slashbox}
\usepackage{subfigure}
\usepackage{epstopdf}
\usepackage{pdfpages}
\usepackage{algorithm}
\usepackage{algorithmic}
\usepackage{url}
\newcommand{\thmlist}{
\begin{list}{Step 1}
{\setlength{\leftmargin}{0.6 in}\setlength{\labelwidth} {0.5 in}}}

\newcommand{\alglist}{
\begin{list}{Step 1}
{\setlength{\leftmargin}{1.1 in} \setlength{\labelwidth}{1.0 in}}}

 \renewcommand{\proof} {\noindent {\bf Proof.} \quad}

\renewcommand{\subtitle}[1]{\color{blue}}

\begin{document}


\title{Continuation methods with the trusty
time-stepping scheme for linearly constrained optimization with noisy data}
\titlerunning{Continuation methods for Linearly Constrained Optimization}
\author{Xin-long Luo\textsuperscript{$\ast$} \and Jia-hui Lv \and Geng Sun}
\authorrunning{Luo, Lv and Sun}

\institute{Xin-long Luo
             \at
              Corresponding author. School of Artificial Intelligence, \\
Beijing University of Posts and Telecommunications, P. O. Box 101, \\
Xitucheng Road  No. 10, Haidian District, 100876, Beijing China\\
             \email{luoxinlong@bupt.edu.cn}            
        \and
        Jia-hui Lv \at
        School of Artificial Intelligence, \\
Beijing University of Posts and Telecommunications, P. O. Box 101, \\
Xitucheng Road  No. 10, Haidian District, 100876, Beijing China\\
             \email{jhlv@bupt.edu.cn}
           \and
           Geng Sun \at
              Institute of Mathematics,
Academy of Mathematics and Systems Science, \\
Chinese Academy of Sciences, 100190, Beijing China\\
    \email{sung@amss.ac.cn}
}

\date{Received: date / Accepted: date}
\maketitle

\begin{abstract}
  The nonlinear optimization problem with linear constraints has many
  applications in engineering fields such as the visual-inertial navigation
  and localization of an unmanned aerial vehicle maintaining the horizontal
  flight. In order to solve this practical problem efficiently, this paper
  constructs a continuation method with the trusty time-stepping scheme
  for the linearly equality-constrained optimization problem at every sampling
  time. At every iteration, the new method only solves a system of linear equations
  other than the traditional optimization method such as the sequential quadratic
  programming  (SQP) method, which needs to solve a quadratic programming
  subproblem. Consequently, the new method can save much more computational
  time than SQP. Numerical results show that the new method works well for this problem
  and its consumed time is about one fifth of that of SQP (the built-in subroutine
  fmincon.m of the MATLAB2018a environment) or that of the traditional dynamical
  method (the built-in subroutine ode15s.m of the MATLAB2018a environment).
  Furthermore, we also give the global convergence analysis of the new method.
\end{abstract}



\vskip 2mm




\section{Introduction} \label{SUBINT}

\vskip 2mm

In this article, we consider the following linearly  equality-constrained
optimization problem
\begin{align}
  &\min_{x \in \Re^n} \; f(x)  \nonumber \\
  &\text{subject to} \; \; Ax = b,   \label{LEQOPT}
\end{align}
where matrix $A \in \Re^{m \times n}$ and vector $b \in \Re^{m}$
may have random noise. This problem has many applications in engineering
fields such as the visual-inertial navigation of an unmanned aerial
vehicle maintaining the horizontal flight \cite{CMFO2009,LLS2020,ZS2015},
and there are many practical methods to solve it such as the sequential
quadratic programming (SQP) method \cite{NW1999} or the penalty function
method \cite{FM1990}.

\vskip 2mm

The penalty function method obtains the solution of the linearly equality-constrained
optimization problem \eqref{LEQOPT} via solving  the following sequential unconstrained
minimization
\begin{align}
   \min_{x\in \Re^{n}} \; P_{\sigma}(x) = f(x) + \sigma \|Ax-b\|^2,
   \label{PENFX}
\end{align}
with increasing $\sigma$. If we denote the global optimal solution of the
unconstrained optimization problem \eqref{PENFX} as $x_{\sigma}^{\ast}$, it is well known
that
\begin{align}
  \lim_{\sigma \to \infty} \; x_{\sigma}^{\ast} = x^{\ast}, \nonumber
\end{align}
where $x^{\ast}$ is the optimal solution of the original constrained optimization
problem \eqref{LEQOPT} \cite{FM1990}. The penalty function method has the asymptotic
convergence as $\sigma \to \infty$ for the constrained optimization problem \eqref{LEQOPT}.
However, in practice, it will meet the ill conditioning which depends on the ratio
of the largest to the smallest eigenvalue (the condition number) of the Hessian matrix
$\nabla_{xx}^{2} P_{\sigma}(x_{\sigma}^{\ast})$, and this ratio tends to increase with
$\sigma$ (pp. 475-476, \cite{Bertsekas2018}). It can be roughly shown as follows.

\vskip 2mm

We denote the rank of matrix $A$ as $r$ and assume that $r < \min\{m, \, n\}$.
From problem \eqref{PENFX}, we obtain the Hessian matrix $ H_{\sigma}(x)$ of
$P_{\sigma}(x)$ via the simple calculation as follows:
\begin{align}
   H_{\sigma}(x) = \nabla_{xx}^{2} P_{\sigma}(x) = \nabla^{2}f(x)
    + 2\sigma A^{T}A.    \label{HESSMATP}
\end{align}
We define $\mu_{i}(B) \, (i = 1, \, 2, \ldots, \, n)$ as the eigenvalues of
matrix $B \in \Re^{n \times n}$. $\mu_{min}(B)$ and $\mu_{max}(B)$ represent
the smallest and largest eigenvalues of matrix $B$, respectively.
From the Courant-Fisher minimax theorem (p. 441, \cite{GV2013}) and
equation \eqref{HESSMATP}, we have
\begin{align}
    & \mu_{min} \left(\nabla^{2}f(x)\right) \le  \mu_{min}(H_{\sigma}(x))  = \min_{\|y\|=1}y^{T}H_{\sigma}(x) y
    \le \min_{\|y\|=1, \; Ay = 0} y^{T}H_{\sigma}(x) y \nonumber \\
    & \quad = \min_{\|y\| =1, \; Ay =0} y^{T}\nabla^{2}f(x) y \le
    \mu_{max}\left(\nabla^{2}f(x)\right)
    \le \max_{i= 1, \, 2, \, \ldots, \, n}\left|\mu_{i}\left(\nabla^{2}f(x)\right)\right|
    \triangleq M(x). \label{MINEIGHLEMX}
\end{align}
By combining $\mu_{min} \left(\nabla^{2}f(x)\right) \ge - M(x)$ with equation 
\eqref{MINEIGHLEMX}, we have
\begin{align}
  \min_{1\le i \le n} \left|\mu_{i}(H_{\sigma}(x))\right|
  \le \left|\mu_{min}(H_{\sigma}(x))\right| \le M(x). \label{MINEIGHX}
\end{align}

\vskip 2mm

Similarly, from equation \eqref{HESSMATP}, we have
\begin{align}
    & \mu_{max}(H_{\sigma}(x)) = \max_{\|y\|=1} y^{T}H_{\sigma}(x)y
     = \max_{\|y\|=1}\left(\sigma y^{T}A^{T}Ay+ y^{T}\nabla^{2}f(x)y \right)
     \nonumber \\
    & \quad \ge \max_{\|y\|=1}\left(\sigma y^{T}A^{T}Ay +
     \min_{\|y\|=1}y^{T}\nabla^{2}f(x)y \right) = \max_{\|y\|=1} \sigma y^{T}A^{T}Ay
    + \min_{\|y\|=1}y^{T}\nabla^{2}f(x)y
    \nonumber \\
    & \quad = \sigma \mu_{max}\left(A^{T}A\right) + \min_{\|y\|=1}y^{T}\nabla^{2}f(x)y
    = \sigma \mu_{max}\left(A^{T}A\right) + \mu_{min}\left(\nabla^{2}f(x)\right).
  \label{MAXEIGHX}
\end{align}
From equations \eqref{MINEIGHX}-\eqref{MAXEIGHX}, we obtain
\begin{align}
  \lim_{\sigma \to \infty}
  \frac  {\mu_{max}(H_{\sigma}(x))}
  {\min_{1 \le i \le n}|\mu_{i}(H_{\sigma}(x))|} \ge \lim_{\sigma \to \infty}
  \frac{\sigma \mu_{max}(A^{T}A) + \mu_{min}\left(\nabla^{2}f(x)\right)}
  {M(x)}  = \infty. \label{ILLCONHX}
\end{align}
That is to say, the condition number of the Hessian matrix $H_{\sigma}(x)$ tends
to infinity.

\vskip 2mm

In order to overcome the numerical difficulty of the penalty function method
near the optimal point $x^{\ast}$ of the constrained optimization problem
\eqref{LEQOPT}, there are some promising methods such as the dynamical methods \cite{AG2003,CKK2003,Goh2011,KLQCRW2008,Tanabe1980} or the SQP methods
\cite{Bertsekas2018,Heinken1996,NW1999} for this problem via handling its first-order
Karush-Kuhn-Tucker conditions directly. The advantage of the dynamical method over
the SQP method is that the dynamical method is capable of finding many local optimal
points of non-convex optimization problems by tracking the trajectories, and it is
even possible to find the global optimal solution \cite{BB1989,Schropp2000,Yamashita1980}.
However, the dynamical method requires more iteration steps and consumes more time
than SQP. In order to improve the computational efficiency of the dynamical method, 
we consider a continuation method with the new time-stepping scheme based on the 
trust-region technique in this article.

\vskip 2mm

The rest of the paper is organized as follows. In section 2, we construct
a new continuation method with the trusty time-stepping scheme
for the linearly equality-constrained optimization problem \eqref{LEQOPT}.
In section 3, we give the global convergence analysis of this new method.
In section 4, we report some promising numerical results of the new method,
in comparison to the traditional SQP method and the traditional dynamical
method for some large scale test problems and a real-world optimization
problem which arises from the visual-inertial navigation and localization problem
with or without the random errors. Finally, we give some discussions and
conclusions in section 5.

\vskip 2mm

\section{Continuation Methods with the Trusty Time-stepping Scheme}

\vskip 2mm

In this section, we construct a continuation method with the new
time-stepping scheme based on the trust-region technique \cite{Yuan2015}
for the linearly equality-constrained optimization problem \eqref{LEQOPT}
via following the trajectory of the differential-algebraic dynamical system
to obtain its equilibrium point.

\vskip 2mm

\subsection{The Differential-Algebraic Dynamical System}

\vskip 2mm

For the linearly constrained optimization problem \eqref{LEQOPT}, it is well known
that its optimal solution $x^{\ast}$ needs to satisfy the Karush-Kuhn-Tucker
conditions (p. 328, \cite{NW1999}) as follows:
\begin{align}
  \nabla_{x} L(x, \lambda) &= \nabla f(x) + A^{T} \lambda = 0,
    \label{FOKKTG} \\
  Ax - b & = 0,             \label{FOKKTC}
\end{align}
where the Lagrangian function $L(x, \lambda)$ is defined by
\begin{align}
  L(x, \; \lambda) = f(x) + \lambda^{T}(Ax-b).
      \label{LAGFUN}
\end{align}
Similarly to the method of the negative gradient flow for the unconstrained
optimization problem \cite{LKLT2009}, from the first-order necessary conditions
\eqref{FOKKTG}-\eqref{FOKKTC}, we can construct a dynamical system of
differential-algebraic equations for problem \eqref{LEQOPT}
\cite{LL2010,Luo2012,LLW2013,Schropp2003} as follows:
\begin{align}
    & \frac{dx}{dt} = - \nabla L_{x}(x, \lambda)
      = -\left(\nabla f(x) + A^{T} \lambda \right),  \label{DAGF} \\
    & Ax - b = 0.                      \label{LACON}
\end{align}

\vskip 2mm

By differentiating the algebraic constraint \eqref{LACON} with respect to $t$
and replacing it into the differential equation \eqref{DAGF}, we obtain
\begin{align}
  A\frac{dx}{dt} = - A \left(\nabla f(x) + A^{T} \lambda \right)
  = - A \nabla f(x) - AA^{T} \lambda = 0.    \label{DIFALGC}
\end{align}
If we assume that matrix $A$ has full row rank further, from equation
\eqref{DIFALGC}, we obtain
\begin{align}
   \lambda = - \left(AA^{T} \right)^{-1} A \nabla f(x). \label{LAMBDA}
\end{align}
By replacing $\lambda$ of equation \eqref{LAMBDA} into equation \eqref{DAGF},
we obtain the dynamical system of ordinary differential equations (ODEs) as
follows:
\begin{align}
  \frac{dx}{dt} = - \left( I - A^{T} \left(AA^{T}\right)^{-1}A\right)
  \nabla f(x). \label{ODGF}
\end{align}
Thus, we also obtain the projection gradient flow for the constrained optimization
problem \cite{Tanabe1980}.

\vskip 2mm

For convenience, we denote the projection matrix $P$ as
\begin{align}
  P  = I - A^{T} \left(AA^{T}\right)^{-1}A.  \label{PROMAT}
\end{align}
It is not difficult to verify $P^{2} = P$ and $PA^{T} = 0$. That is to say, $P$
is a symmetric projection matrix and its eigenvalues are 0 or 1. From Theorem 2.3.1
in p. 73 of \cite{GV2013}, we know that its matrix 2-norm is
\begin{align}
     \|P\| = 1. \label{MATNP}
\end{align}

\vskip 2mm

\begin{remark}
If $x(t)$ is the solution of ODEs \eqref{ODGF}, it is not difficult
to verify that $x(t)$ satisfies $A (dx/dt) = 0$. That is to say, if the initial point
$x_{0}$ of ODEs \eqref{ODGF} satisfies $Ax_{0} = b$, the solution $x(t)$ of ODEs
\eqref{ODGF} also satisfies $Ax(t) = b, \; \forall t \ge 0$.
\end{remark}

\vskip 2mm 

\begin{remark}
If we assume that $x(t)$ is the solution of ODEs \eqref{ODGF}, from equations
\eqref{PROMAT}-\eqref{MATNP}, we obtain
\begin{align}
 \frac{df(x)}{dt} & = \left(\nabla f(x)\right)^{T} \frac{dx}{dt}
 = - (\nabla f(x))^{T} P \nabla f(x) = -
 (\nabla f(x))^{T} P^{2} \nabla f(x) \nonumber \\
 & = - (P \nabla f(x))^{T}(P \nabla f(x)) = - \|P \nabla f(x) \|_{2}^{2} \le 0.
 \nonumber
\end{align}
That is to say, $f(x)$ is monotonically decreasing along the solution curve $x(t)$
of the dynamical system \eqref{ODGF}. Furthermore,
the solution $x(t)$ converges to $x^{\ast}$ when $t$ tends to infinity
\cite{Schropp2000,Tanabe1980}, where $x^{\ast}$ satisfies the first-order
Karush-Kuhn-Tucker conditions \eqref{FOKKTG}-\eqref{FOKKTC}. Thus, we can follow
the trajectory $x(t)$ of the ordinary differential equation \eqref{ODGF} or
the trajectory $(x(t), \, \lambda(t))$ of differential-algebraic equations
\eqref{DAGF}-\eqref{LACON} to obtain their equilibrium point $x^{\ast}$, which
is also one saddle point of the original optimization problem \eqref{LEQOPT}.
\end{remark}

\vskip 2mm

\subsection{Continuation Methods} \label{SUBSICM}

\vskip 2mm

The solution curve of general differential-algebraic equations is not efficiently
followed on an infinite interval by the traditional ODE method
\cite{AP1998,BCP1996,HW1996,LF2000}, so one needs to construct the
particular method for this problem \eqref{DAGF}-\eqref{LACON}. We regard the
algebraic equation \eqref{LACON} as a degenerate differential equation
\cite{BCP1996,CKK2003,Shampine2002}, and apply the first-order implicit Euler
method to the system of differential-algebraic equations \eqref{DAGF}-\eqref{LACON},
then we obtain
\begin{align}
   & x_{k+1} =  x_{k} - \Delta t_{k} \left( \nabla f(x_{k+1})
   + A^{T}\lambda_{k+1} \right),      \label{IEDD} \\
   & Ax_{k+1} - b = 0,      \label{IEDA}
\end{align}
where $\Delta t_k$ is the time-stepping size.

\vskip 2mm

Since the system of equations \eqref{IEDD}-\eqref{IEDA} is a nonlinear system
which is not directly solved, we seek for its explicit approximation formula.
We denote $s_{k} = x_{k+1} - x_{k}$. By using the first-order Taylor expansion,
we have the linear approximation $\nabla f(x_{k}) + \nabla^{2} f(x_{k})s_{k}$ of
$\nabla f(x_{k+1})$. From equation \eqref{LAMBDA} and the first-order Taylor expansion,
we have
\begin{align}
    \lambda_{k+1} & = - (AA^{T})^{-1}A\nabla f(x_{k+1})
    \approx -(AA^{T})^{-1}A \nabla f(x_{k}) - (AA^{T})^{-1}A \nabla^{2} f(x_{k})s_{k}
    \nonumber \\
    & = \lambda_{k} - (AA^{T})^{-1}A \nabla^{2} f(x_{k})s_{k}. \nonumber
\end{align}
By replacing them into equation \eqref{IEDD}, we obtain the predictor $x_{k+1}^{P}$
of $x_{k+1}$ as follows:
\begin{align}
    \left({1}/{\Delta t_{k}}I + B_{k} \right)d_{k}
    & = - p_{g_{k}}, \label{PRDK} \\
    x_{k+1}^{P} & = x_{k} + d_{k}, \label{PRXK1}
\end{align}
where $B_{k}$ equals the Jacobian matrix $\nabla^{2} f(x_{k}) + A^{T}
\frac{\partial}{\partial x}\lambda(x_{k})$ ($\lambda(x)$ defined by equation
\eqref{LAMBDA}) or its quasi-Newton approximation matrix, and
\begin{align}
    p_{g_{k}} = \nabla_{x} L(x_k, \; \lambda_k)
    = \nabla f(x_{k})+ A^{T}\lambda_{k}. \label{GRADLGK}
\end{align}

\vskip 2mm

The predicted point $x_{k+1}^{P}$ will escape from the constraint plane \eqref{LACON},
so we pull it back to the constraint plane by solving the following projection problem:
\begin{align}
     \min_{x \in \Re^{n}} \; \left\|x - x_{k+1}^{P}\right\|^{2}
     \; \text{subject to} \hskip 2mm  Ax = b.  \label{MINDIST}
\end{align}
It is not difficult to obtain the solution of the linearly constrained least-squares
problem \eqref{MINDIST} via using the Lagrangian multiplier method
(p. 479, \cite{Bertsekas2018}) as follows:
\begin{align}
     x_{k+1} = x_{k+1}^{P} + A^{T}\left(AA^{T}\right)^{-1} \left(b - Ax_{k+1}^{P}\right).
     \label{ITK1XK1}
\end{align}
Notice that $x_{k+1}$ satisfies the constraint $Ax = b$. From equation \eqref{PRXK1}
and equation \eqref{ITK1XK1}, we have
\begin{align}
    x_{k+1} & = x_{k+1}^{P} + A^{T}\left(AA^{T}\right)^{-1} \left(b - Ax_{k+1}^{P}\right)
       = x_{k} + d_{k} + A^{T}\left(AA^{T}\right)^{-1} \left(Ax_{k} - Ax_{k+1}^{P}\right)
      \nonumber \\
    & = x_{k} + d_{k} - A^{T}\left(AA^{T}\right)^{-1} Ad_{k}
      = x_{k} + \left(I - A^{T}\left(AA^{T}\right)^{-1} A\right)d_{k} \nonumber \\
    & = x_{k} + Pd_{k}, \label{PRODK}
\end{align}
where the projection matrix $P$ is defined by equation \eqref{PROMAT}.

\vskip 2mm

After solving $x_{k+1}$ from equation \eqref{PRDK} and equation \eqref{PRODK},
according to equation \eqref{LAMBDA}, we obtain the Lagrangian multiplier
$\lambda_{k+1}$ as follows:
\begin{align}
     \lambda_{k+1} = - \left(AA^{T}\right)^{-1} A\nabla f(x_{k+1}).  \label{LAMBDAK1}
\end{align}
By replacing $\lambda_{k+1}$ of equation \eqref{LAMBDAK1} into equation \eqref{GRADLGK},
we also obtain
\begin{align}
    p_{g_k} = \nabla_{x} L(x_k, \; \lambda_k)
    = \nabla f(x_{k})+ A^{T}\lambda_{k} =
    \left(I -A^{T}\left(AA^{T}\right)^{-1}A\right) \nabla f(x_{k}) = Pg_{k},
    \label{PROGK}
\end{align}
where $g_{k} = \nabla f(x_{k})$ and the projection matrix $P$ is defined by equation
\eqref{PROMAT}.

\vskip 2mm

\subsection{The Trusty Time-stepping Scheme}

\vskip 2mm

Another issue is how to adaptively adjust the time-stepping size $\Delta t_k$
at every iteration. We borrow the adjustment method of the trust-region radius
from the trust-region method due to its robust convergence and fast local
convergence \cite{CGT2000}. Since $x_{k+1}$ is the solution of the linearly
constrained least-squares problem \eqref{MINDIST}, $x_{k+1}$ maintains the feasibility.
Therefore, we use the objective function $f(x)$ instead of the nonsmooth penalty
function $f(x) + \sigma \|Ax-b\|_{1}$ as the cost function.

\vskip 2mm

When we use the trust-region technique to adaptively adjust time-stepping size
$\Delta t_{k}$ \cite{Higham1999}, we also need to construct a local approximation
model of the objective $f(x)$ around $x_{k}$. Here, we adopt the following quadratic
function as its approximation model:
\begin{align}
  q_k(x) = f(x_k) + (x-x_{k})^{T} g_{k} + {1}/{2}(x-x_k)^{T}
      B_{k} (x-x_{k}). \label{SOAM}
\end{align}
where $g_{k} = \nabla f(x_k)$ and the symmetric matrix $B_{k}$ equals
$\nabla^{2} f(x_{k}) + A^{T} \frac{\partial}{\partial x} \lambda(x_{k})$ ($\lambda$ defined by equation
\eqref{LAMBDA}) or its quasi-Newton approximation matrix.
We enlarge or reduce the time-stepping size $\Delta t_k$ at every iteration
according to the following ratio:
\begin{align}
    \rho_k = \frac{f(x_k)-f(x_{k+1})}{q_k(x_k)-q_k(x_{k+1})}.
    \label{MRHOK}
\end{align}
A particular adjustment strategy is given as follows:
\begin{align}
     \Delta t_{k+1} = \begin{cases}
          \gamma_1 \Delta t_k, &{if \hskip 1mm 0 \leq \left|1- \rho_k \right| \le \eta_1,}\\
          \Delta t_k, &{if \hskip 1mm \eta_1 < \left|1 - \rho_k \right| < \eta_2,}\\
          \gamma_2 \Delta t_k, &{if \hskip 1mm \left|1-\rho_k \right| \geq \eta_2,}
                   \end{cases} \label{ADTK1}
\end{align}
where the constants are selected as $\eta_1 = 0.25, \; \gamma_1 = 2, \; \eta_2 = 0.75, \;
\gamma_2 = 0.5$  according to numerical experiments.

\vskip 2mm

\subsection{The Treatments of Deficient Rank and Infeasible Initial Points}
\label{SUBSECDEG}

\vskip 2mm

For a real-world problem, the rank of matrix $A$ may be deficient and even
the constraint system may be inconsistent when the data ($A,\, b$) have random noise.
We handle this problem via solving the following best approximation problem
\begin{align}
    \min_{x \in \Re^{n}} \; \|Ax - b\|^{2} \label{LLSP}
\end{align}
to obtain the reduced constraint system of problem \eqref{LEQOPT}.

\vskip 2mm

First, we factorize matrix $A$ with its singular value decomposition
(pp. 76-80, \cite{GV2013}) as follows:
\begin{align}
   A = U \Sigma V^{T}, \; \Sigma = diag(\sigma_1, \, \sigma_2, \, \ldots,
   \, \sigma_r, 0, \, \ldots, \, 0),  \;
   \sigma_1 \ge \sigma_2 \ge \ldots \ge \sigma_r > 0, \label{SVDOFA}
\end{align}
where $U \in \Re^{m \times m}$ and $V \in \Re^{n \times n}$ are orthogonal
matrices, and $r$ is the rank of matrix $A$. Thus, problem \eqref{LLSP} equals the
following linear least-squares problem
\begin{align}
  \min_{x \in \Re^{n}} \; \| \Sigma V^{T}x - U^{T}b\|^{2}, \label{APPCON}
\end{align}
which leads to the reduced constraint system
\begin{align}
   V_{r}^{T} x = b_{r},  \label{APPCON}
\end{align}
where $V_{r} = V(1:n, 1:r)$, \, $V_{r}^{T}V_{r} = I$,
and $b_{r} = ((U^{T}b)(1:r))./diag(\Sigma(1:r, \, 1:r))$.

\vskip 2mm

Therefore, when the constraint system of problem \eqref{LEQOPT} is
consistent, it equals the following problem
\begin{align}
   \min_{x \in \Re^{n}} \; f(x) \; \text{subject to} \; V_{r}^{T} x = b_{r}.
   \label{NLOBJFUN}
\end{align}
When the constraint system of problem \eqref{LEQOPT} is inconsistent,
problem \eqref{NLOBJFUN} is the best relaxation approximation
of the original optimization problem \eqref{LEQOPT}. After this preprocess,
we reformulate the projection matrix $P$ defined by equation \eqref{PROMAT}
as follows:
\begin{align}
   P = I - V_{r}V_{r}^{T}. \label{SMPPROJ}
\end{align}

\vskip 2mm

Consequently, in subsection \ref{SUBSICM}, we only need to replace matrix $A$
and vector $b$ with matrix $V_{r}^{T}$ and vector $b_{r}$ respectively, then the
continuation method \eqref{PRODK} can handle the deficient rank problem.

\vskip 2mm

For a real-world optimization problem \eqref{LEQOPT}, we probably meet the infeasible
initial point $x_{0}$. That is to say, the initial point can not satisfy the constraint
\eqref{APPCON}. We handle this problem by solving the following projection problem:
\begin{align}
     \min_{x \in \Re^{n}} \; \left\|x - x_{0} \right\|^{2}
     \; \text{subject to} \hskip 2mm  V_{r}^{T} x = b_{r},  \label{MINDISTVB}
\end{align}
where $V_{r} \in \Re^{n \times r}$ satisfies $V_{r}^{T}V_{r} = I$. By using the
Lagrangian multiplier method to solve problem \eqref{MINDISTVB}, we obtain the
initial feasible point $x_{0}^{F}$ of problem \eqref{NLOBJFUN} as follows:
\begin{align}
    x_{0}^{F} = x_{0} + V_{r}\left(b_{r} - V_{r}^{T} x_{0}\right).
    \label{INFIPT}
\end{align}
For convenience, we set $x_{0} = x_{0}^{F}$ in line 4, Algorithm \ref{ALGPTCTR}.

\vskip 2mm

\subsection{The BFGS Quasi-Newton Updating Method}

\vskip 2mm

For the large-scale problem, the numerical estimation of the Hessian matrix
$\nabla^{2}f(x_{k})$ consumes much time. In order to overcome this shortcoming,
we use the BFGS quasi-Newton formula (pp. 194-198, \cite{NW1999})
to update the approximation $B_{k}$ of $\nabla^{2} f(x_{k}) + A^{T} \frac{\partial}{\partial x} \lambda(x_{k})$,
where $\lambda(x)$ is defined by equation \eqref{LAMBDA}. The BFGS updating formula
can be written as
\begin{align}
  B_{k+1} = B_{k} - \frac{B_{k}s_{k}s_{k}^{T}B_{k}}{s_{k}^{T}B_{k}s_{k}}
  + \frac{y_{k}y_{k}^{T}}{y_{k}^{T}s_{k}}, \label{BFGS}
\end{align}
where $s_{k} = x_{k+1} - x_{k}$ and $y_{k} = \nabla_{x} L(x_{k+1}, \, \lambda_{k+1})
- \nabla_{x} L(x_{k}, \, \lambda_{k})$ ($\nabla_{x}L(x_{k}, \, \lambda_{k})$ defined
by equation \eqref{GRADLGK}). The initial matrix $B_{0}$ can be simply selected
by the identity matrix.

\vskip 2mm

The BFGS updating formula \eqref{BFGS} has some nice properties such as the
symmetric positive definite property of matrix $B_{k+1}$ if $B_{k}$ is symmetric
positive definite and $y_{k}^{T}s_{k} > 0$. Its proof can be found in p. 199,
\cite{NW1999}. For convenience, we state its brief proof as follows.

\vskip 2mm

Since matrix $B_{k}$ is symmetric positive definite, we have its Cholesky
factorization $B_{k} = L_{k}L_{k}^{T}$ and denote
\begin{align}
    \alpha_{k} = L_{k}^{T}s_{k}, \; \text{and} \; \beta_{k} = L_{k}^{T}z
    \; \text{for any} \; z \in \Re^{n}. \nonumber
\end{align}
Thus, for any nonzero $z \in \Re^{n}$, from equation \eqref{BFGS} and the assumption
$y_{k}^{T}s_{k} > 0$, we have
\begin{align}
    z^{T}B_{k+1}z  = z^{T}B_{k}z - \frac{(z^{T}B_{k}s_{k})^{2}}{s_{k}^{T}B_{k}s_{k}}
     + \frac{(z^{T}y_{k})^{2}}{y_{k}^{T}s_{k}}
     = \|\beta_{k}\|^{2} - \frac{(\alpha_{k}^{T}\beta_{k})^{2}}{\|\alpha_{k}\|^{2}}
     + \frac{(z^{T}y_{k})^{2}}{y_{k}^{T}s_{k}} \ge 0.  \label{NONND}
\end{align}
In the last inequality of equation \eqref{NONND}, we use the Cauchy-Schwartz
inequality $\|\alpha_{k}\| \|\beta_{k} \|$ $\ge |\alpha_{k}^{T}\beta_{k}|$.
Its equality holds if only if $\beta_{k} =  t \alpha_{k}$. In this case, we have
$s_{k} = t z$ since $\alpha_{k} = L_{k}^{T}s_{k} = t \beta_{k} = t L_{k}^{T}z$.
Thus, we have $(z^{T}y_{k})^{2}/y_{k}^{T}s_{k} = t^{2} > 0$, which leads to
$z^{T}B_{k+1}z > 0$.

\vskip 2mm

According to the above discussions, we give the detailed implementation of
the continuation method with the trusty time-stepping scheme for the
linearly equality-constrained optimization problem \eqref{LEQOPT} in Algorithm
\ref{ALGPTCTR}.

\begin{algorithm}
	\renewcommand{\algorithmicrequire}{\textbf{Input:}}
	\renewcommand{\algorithmicensure}{\textbf{Output:}}
    \newcommand{\algorithmicbreak}{\textbf{break}}
    \newcommand{\BREAK}{\STATE \algorithmicbreak}
	\caption{Continuation methods with the trusty time-stepping scheme for
            the linearly equality-constrained optimization problem (the Ptctr method)}
    \label{ALGPTCTR}	
	\begin{algorithmic}[1]
		\REQUIRE ~~\\
        the objective function: $f(x)$; \\
        the linear constraint: $Ax  = b$; \\
        the initial point: $x_0$ \, (optional); \\
        the terminated parameter: $\epsilon$ \, (optional).
		\ENSURE ~~\\
        the optimal approximation solution $x^{\ast}$.

        \vskip 2mm
        		
        \STATE If the called function does not provide the initial values $x_0$ and
        $\epsilon$, we set $x_0 = [1,\, 1, \, \ldots, \, 1]^{T}$ and
        $\epsilon = 10^{-6}$, respectively.
        \STATE Initialize the parameters: $\eta_{a} = 10^{-6}, \; \eta_1 = 0.25,
        \; \gamma_1 =2, \; \eta_2 = 0.75, \; \gamma_2 = 0.5$.
        \STATE Factorize matrix $A$ with the singular value decomposition as follows:
        $$
         A = U \Sigma V^{T}, \; \Sigma = diag(\sigma_1, \, \sigma_2, \, \ldots,
          \, \sigma_r, \, 0, \, \ldots, \, 0),
        $$
        and denote $V_{r} = V(1:n, 1:r)$,
        $b_{r} = ((U^{T}b)(1:r))./diag(\Sigma(1:r, \, 1:r))$, where $r$ is the rank
        of matrix $A$.
        \STATE Compute
        $$ x_{0} \leftarrow x_{0} + V_{r}\left(b_{r} - V_{r}^{T} x_{0}\right),
        $$
        such that $x_{0}$ satisfies the linear system of constraints $V_{r}^{T}x = b_{r}$.
        \STATE Set $k = 0$ and $B_{0} = I$. Evaluate $f_0 = f(x_0)$ and $g_0 = \nabla f(x_0)$.
        \STATE Compute the projection gradient $p_{g_{0}} = g_{0} - V_{r}\left(V_{r}^{T}g_{0}\right)$.
        \STATE Compute the initial time-stepping size
        $$ \Delta t_0 = \min\left\{10^{-2}, \; {1}/{\|p_{g_0}\|}\right\}. $$
        \WHILE{$\|p_{g_k}\|> \epsilon$}
          \IF{$\left(\frac{1}{\Delta t_{k}}I + B_k - P^{T}B_{k}P\right) \succ 0$ \&\&
          $\left(\frac{1}{\Delta t_{k}}I + B_k\right) \succ 0$}
            \STATE Solve the linear system \eqref{PRDK} by the Cholesky factorization
            to obtain the search direction $d_{k}$, and compute
            $$
              x_{k+1} = x_{k} + Pd_{k} = x_{k} + \left(I - V_{r}V_{r}^{T}\right)d_{k}
              = x_{k} + d_{k} - V_{r}\left(V_{r}^{T}d_{k}\right).
            $$
            \STATE Evaluate $f_{k+1} = f(x_{k+1})$ and compute the ratio $\rho_{k}$
            from equations \eqref{SOAM}-\eqref{MRHOK}.
          \ELSE
            \STATE Let $\rho_k  = -1$.
          \ENDIF
          \IF{$\rho_k\le \eta_{a}$}
            \STATE Set $x_{k+1} = x_{k}, \; f_{k+1} = f_{k}, \; p_{g_{k+1}} = p_{g_{k}},
            \; g_{k+1} = g_{k}, \; B_{k+1} = B_{k}.$
          \ELSE
            \STATE Evaluate $g_{k+1} = \nabla f(x_{k+1})$, and
            update $B_{k+1}$ by the BFGS formula \eqref{BFGS}.
            \STATE Compute the projection gradient
            $$
              p_{g_{k+1}} = Pg_{k+1} = \left(I - V_{r}V_{r}^{T}\right)g_{k+1}
              = g_{k+1} - V_{r}\left(V_{r}^{T}g_{k+1}\right).
            $$
          \ENDIF
          \STATE Adjust the time-stepping size $\Delta t_{k+1}$ based on the
          trust-region updating scheme \eqref{ADTK1}.
          \STATE Set $k \leftarrow k+1$.
        \ENDWHILE
	\end{algorithmic}
\end{algorithm}

\section{Algorithm Analysis}

In this section, we analyze the global convergence of the continuation method
with the trusty time-stepping scheme for the linearly equality-constrained
optimization problem (i.e. Algorithm \ref{ALGPTCTR}). Firstly, we give a lower-bounded
estimate of $q_{k}(x_k) - q_{k}(x_{k+1})$ $(k = 1, \, 2, \, \ldots)$. This
result is similar to that of the trust-region method for the unconstrained
optimization problem \cite{Powell1975}. For simplicity, we assume that the rank of
matrix $A$ is full and the constraint $Ax = b$ is consistent.

\vskip 2mm

\begin{lemma} \label{LBSOAM}
Assume that the quadratic model $q_{k}(x)$ is defined by equation \eqref{SOAM}
and $d_{k}$ is the solution of equation \eqref{PRDK}. Furthermore, we suppose that
the time-stepping size $\Delta t_k$ satisfies
\begin{align}
    \left(1/{\Delta t_{k}} \, I + B_k\right) \succ 0 \;
    \text{and} \; \left({1}/{\Delta t_{k}} \, I + B_{k}- P^{T}B_{k}P \right) \succ 0,
    \label{PSDASS}
\end{align}
where matrix $P$ is defined by equation \eqref{SMPPROJ}.
Then, we have
\begin{align}
    q_{k}(x_k) - q_{k}(x_{k}+Pd_{k}) \ge {1}/{2} \left\|p_{g_{k}} \right\|
    \min \left\{\left\|Pd_{k}\right\|, \; \|p_{g_k}\| /(3\|B_{k}\|)\right\},
    \label{PLBREDST}
\end{align}
where the projection gradient $p_{g_k} = Pg_{k}$.
\end{lemma}
\proof Let $\tau_{k} = 1/\Delta t_{k}$. From equation \eqref{PRDK}, we obtain
\begin{align}
    q_{k}(x_{k}) &- q_{k}(x_{k}+Pd_{k}) = - {1}/{2}\left(d_{k}^{T}P^{T}B_{k}Pd_{k}\right)
     - (Pg_{k})^{T}d_{k}  \nonumber \\
    & = - {1}/{2}\left(d_{k}^{T}P^{T}B_{k}Pd_{k}\right)
    + p_{g_k}^{T} \left(\tau_{k}I + B_{k}\right)^{-1}p_{g_k}  \nonumber \\
    & =  {1}/{2}\left(p_{g_k}^{T} \left(\tau_{k}I + B_{k}\right)^{-1}p_{g_k}
    + d_{k}^{T}\left(-P^{T}B_{k}P + \tau_{k}I + B_{k}\right)d_{k}\right).
    \label{DKGKDK}
\end{align}
We denote $\mu_{min}\left(B_{k}-P^{T}B_{k}P\right)$ as the smallest eigenvalue of
matrix $\left(B_{k}-P^{T}B_{k}P\right)$, and set
\begin{align}
    \tau_{lb} = \min\left\{0, \; \mu_{min} \left(B_{k}-P^{T}B_{k}P\right)\right\}.
    \label{LOWBTAU}
\end{align}
From equations \eqref{PSDASS}, \eqref{DKGKDK}-\eqref{LOWBTAU} and the bound
on the eigenvalues of matrix $(\tau_{k}I + B_{k})^{-1}$, we obtain
\begin{align}
    & q_{k}(x_{k}) - q_{k}(x_{k}+Pd_{k}) \ge {1}/{2}\left(p_{g_k}^{T}
    \left(\tau_{k}I + B_{k}\right)^{-1}p_{g_k}
    + \left(\tau_{k} + \tau_{lb}\right)\left\|d_{k}\right\|^{2}\right)
    \nonumber \\
    & \quad \ge {1}/{2}\left({\left\|p_{g_{k}}\right\|^{2}}/{\left(\tau_{k}+\left\|B_{k}\right\|\right)}
    + \left(\tau_{k} + \tau_{lb}\right)\left\|d_{k}\right\|^{2} \right).
    \label{LBEQK}
\end{align}
In the above second inequality, we use the property $|\mu_{i}(B_{k})|\le \|B_{k}\|$,
where $\mu_{i}(B_{k})$ is an eigenvalue of matrix $B_{k}$.

\vskip 2mm

Now we consider the properties of the function
\begin{align}
   \varphi(\tau) \triangleq \tau \left\|d_k\right\|^{2}
   + {\left\|p_{g_k}\right\|^2}/\left(\tau - \tau_{lb}+\left\|B_k \right\|\right).
   \label{VLF}
\end{align}
It is not difficult to verify that the second-order derivative of $\varphi(\tau)$
is positive when $\left(\tau - \tau_{lb} + \|B_{k}\right\|) > 0$ since
$\varphi^{''}(\tau) = 2\|p_{g_k}\|^2 /\left(\tau - \tau_{lb} + \|B_{k}\|\right)^3  \ge 0$.
Thus, the function $\varphi(\tau)$ attains its minimum $\varphi(\tau_{min})$ when
$\tau_{min}$ satisfies $\varphi^{'}(\tau_{min})=0$ and
$\tau \ge -(-\tau_{lb}+\|B_k\|)$. That is to say, we have
\begin{align}
  \varphi(\tau_{min}) = 2\|p_{g_k}\| \|d_k\|
  + \left(\tau_{lb}-\|B_k\|\right) \|d_k\|^{2},   \label{MINV}
\end{align}
where
\begin{align}
    \tau_{min} = \|p_{g_k}\|/\|d_k\| + \tau_{lb}-\|B_k\|. \label{MINLD}
\end{align}

\vskip 2mm

We prove the property \eqref{PLBREDST} when $\tau_{min} \ge 0$ or
$\tau_{min} < 0$ separately as follows.

\vskip 2mm

(i) When $\left(\|p_{g_k} \|/\|d_k\| + \left(\tau_{lb}-\|B_k\|\right)\right) \ge 0$,
from equation \eqref{MINLD}, we have $\tau_{min} \ge 0$. From the assumption
\eqref{PSDASS} and the definition \eqref{LOWBTAU} of $\tau_{lb}$, we have
$\tau_k \ge -\tau_{lb}$. Thus, from equations \eqref{LBEQK}--\eqref{MINLD},
we obtain
\begin{align}
    & q_k(x_{k})   - q_k(x_{k}+Pd_{k}) \ge {1}/{2}
    \left((\tau_{k}+\tau_{lb}) \|d_k\|^{2} + {\|p_{g_k}\|^2}/\left(\tau_{k} + \|B_k\|\right) \right)
    \nonumber \\
     & \quad =  1/2 \varphi(\tau_{k}+\tau_{lb})
    \ge {1}/{2} \varphi(\tau_{min}) \nonumber \\
    & \quad = {1}/{2} \left(\|p_{g_k}\| \|d_k\|
    + \left(\|p_{g_k}\| \|d_k\|+\left(\tau_{lb}-\|B_k\|\right)\|d_k\|^{2}\right)\right)
   \ge {1}/{2}\|p_{g_k}\| \|d_k\|.   \label{MINQCOV}
\end{align}

\vskip 2mm

(ii) The other case is $\left(\|p_{g_k} \|/\|d_k\| + \left(\tau_{lb}-\|B_k\|\right)\right) < 0$.
In this case, from equation \eqref{MINLD}, we have $\tau_{min} < 0$. It is not
difficult to verify that $\varphi(\tau)$ is monotonically increasing
when $\tau \ge 0$ and $\tau_{min}<0$. From the definition \eqref{LOWBTAU} of $\tau_{lb}$ and
the property \eqref{MATNP}, we have
\begin{align}
    |\tau_{lb}| &\le \left|\mu_{min} \left(B_{k}-P^{T}B_{k}P\right)\right|
    \le \|B_{k} - P^{T}B_{k}P\| \le \|B_{k}\|+\|P^{T}B_{k}P\| \nonumber \\
    & \le \|B_{k}\|+\|P^{T}\|\|B_{k}\|\|P\|= 2\|B_k\|. \nonumber
\end{align}
By using this property and the monotonicity of $\varphi(\tau)$, from equations
\eqref{LBEQK}-\eqref{VLF}, we obtain
\begin{align}
     & q_k(x_{k})  - q_k(x_{k}+Pd_{k})
    \ge {1}/{2}\left((\tau_{k}+\tau_{lb}) \|d_k\|^{2}
    + {\|p_{g_k}\|^2}/\left(\tau_k + \|B_k \|\right)\right) \nonumber \\
    & \quad = \frac{1}{2} \varphi(\tau_{k}+\tau_{lb}) \ge \frac{1}{2} \varphi(0)
    = \frac{1}{2(-\tau_{lb}+\|B_k\|)}\|p_{g_k}\|^2
     \ge \frac{1}{6\|B_k\|}\|p_{g_k}\|^2.   \label{MINQMI}
\end{align}

\vskip 2mm

From equations (\ref{MINQCOV})-(\ref{MINQMI}), we get
\begin{align}
     q_k(x_{k}) - q_k(x_{k}+Pd_{k}) \ge {1}/{2} \|p_{g_k}\|
     \min\left\{\|d_k\|, \;  {\|p_{g_k}\|}/{(3\|B_k\|)}\right\}.
     \label{MINQATR}
\end{align}
By using the property \eqref{MATNP} of matrix $P$, we have
\begin{align}
 \|Pd_{k}\| \le \|P\| \|d_{k}\| = \|d_{k}\|. \label{PROPD}
\end{align}
Therefore, from inequalities \eqref{MINQATR}-\eqref{PROPD}, we obtain the estimate
\eqref{PLBREDST}. \qed

\vskip 2mm

In order to prove that $p_{g_k}$ converges to zero when $k$ tends to infinity,
we need to estimate the lower bound of time-stepping sizes
$\Delta t_{k} \, (k = 1, \, 2, \, \ldots)$ when
$\|p_{g_k}\| \ge \epsilon_{p_g} > 0, \; k = 1, \, 2, \, \ldots$.

\vskip 2mm

\begin{lemma} \label{DTBOUND}
Assume that $f: \; \Re^{n} \rightarrow \Re$ is
twice continuously differentiable and the constrained level set
\begin{align}
     S_{f} = \left\{x: \; f(x) \le f(x_0), \;  Ax = b \right\} \label{LSCONFBD}
\end{align}
is bounded. We assume that the Hessian matrix function $\nabla^{2} f(\cdot)$
is Lipschitz continuous. That is to say, it exists a positive constant $L_{c}$
such that
\begin{align}
     \left\|\nabla^{2} f(x) - \nabla^{2} f(y) \right\| \le L_{c} \|x - y\|, \;
     \forall x, \, y  \in \Re^{n}. \label{LIPSCHCON}
\end{align}
 We suppose that the
sequence $\{x_{k}\}$ is generated by Algorithm \ref{ALGPTCTR} and the quasi-Newton
matrices $B_{k} \, (k=1, \, 2, \, \ldots)$ are bounded. That is to say, it exists a
positive constant $M_{B}$ such that
\begin{align}
    \|B_{k}\| \le M_{B}, \; k = 1, \, 2, \ldots. \label{BOUNDBK}
\end{align}
Furthermore, we assume that it exists a positive constant $\epsilon_{p_g}$ such
that
\begin{align}
    \|p_{g_{k}}\| \ge \epsilon_{p_g} > 0, \; k = 1, \, 2, \ldots, \label{PGKGEPN}
\end{align}
where $p_{g_{k}} = P g_{k}$, \; $g_{k} = \nabla f(x_k)$, and $P$ is defined by
equation \eqref{SMPPROJ}. Then, it exists a positive constant $\delta_{\Delta t}$
such that
\begin{align}
    \Delta t_{k} \ge \gamma_{2} \delta_{\Delta t} > 0, \; k = 1, \; 2, \dots,
    \label{DTGEPN}
\end{align}
where $\Delta t_{k}$ is adaptively adjusted by the trust-region updating scheme
\eqref{SOAM}-\eqref{ADTK1}.
\end{lemma}

\vskip 2mm

\proof Since the level set $S_{f}$ is bounded, according to Proposition A.7 in
pp. 754-755 of reference \cite{Bertsekas2018}, $S_{f}$ is closed. Thus, it exists
two positive constants $M_{g}$ and $M_{G}$ such that
\begin{align}
     \|g_{k}\| \le M_{g}, \;  \|G_{k}\| \le M_{G}, \; k = 1, \; 2, \dots,
     \label{PGGKBD}
\end{align}
where $g_{k} = \nabla f(x_k)$ and $G_{k} = \nabla^{2} f(x_k)$. From equation
\eqref{MATNP}, we know $\|P\| = 1$. By using this property and the assumption
\eqref{BOUNDBK}, we have
\begin{align}
     &\left|\mu_{min}\left(B_{k}-P^{T}B_{k}P\right)\right| \le
     \left\|B_{k}-P^{T}B_{k}P\right\| \nonumber \\
     & \quad \le \|B_{k}\| + \|P^{T}\|\|B_{k}\|\|P\|
     = 2\|B_{k}\| \le 2M_{B}, \; k=1, \, 2,\ldots,  \label{PROGKUPBD}
\end{align}
where $\mu_{min}(B)$ represents the smallest eigenvalue of matrix $B$. Thus,
from equation \eqref{PROGKUPBD}, we obtain
\begin{align}
    & \mu_{min}\left({1}/{\Delta t_{k}}I + B_{k}-P^{T}B_{k}P\right)
      = {1}/{\Delta t_{k}} + \mu_{min}\left(B_{k}-P^{T}B_{k}P\right)
     \nonumber \\
    & \quad \ge {1}/{\Delta t_{k}} - 2M_{B}, \; k=1, \, 2, \ldots.
    \label{SEPROGLBD}
\end{align}
Similarly, from the assumption \eqref{BOUNDBK}, we have
\begin{align}
    \mu_{min} \left({1}/{\Delta t_{k}} I + B_{k}\right)
     =  {1}/{\Delta t_{k}}+ \mu_{min} \left(B_{k}\right)
    \ge {1}/{\Delta t_{k}} - M_{B}, \;  k=1, \, 2, \ldots. \label{SEGLBD}
\end{align}
Therefore, the positive definite conditions of equation \eqref{PSDASS} are
satisfied when $\Delta t_{k} < 1/(2M_{B}) \, (k=1,\,2,\ldots)$.

\vskip 2mm

From a second-order Taylor expansion, we have
\begin{align}
    f(x_{k}+Pd_{k}) =  f(x_{k}) + g_{k}^{T}(Pd_{k})
    + 1/2(Pd_{k})^{T}\nabla^{2}f(\tilde{x}_{k})(Pd_{k}), \label{SOTEFK}
\end{align}
where $\tilde{x}_{k} = x_{k} + \theta_{k} (Pd_{k}), \; 0 \le \theta_{k} \le 1$.
From the Lipschitz continuity \eqref{LIPSCHCON} of $\nabla^{2} f(\cdot)$
and the boundedness \eqref{PGGKBD} of $\nabla^{2} f(x_{k})$, we have
\begin{align}
     \left\|\nabla^{2}f(\tilde{x}_{k})\right\|
     \le \left\|\nabla^{2}f(\tilde{x}_{k}) - \nabla^{2}f(x_{k})\right\|
     + \left\|\nabla^{2}f(x_{k})\right\| \le L_{c} \|Pd_{k}\|+ M_{G}.
     \label{BMSOFK}
\end{align}
Thus, from equations \eqref{MRHOK}, \eqref{PLBREDST}, \eqref{PGGKBD},
\eqref{SOTEFK}-\eqref{BMSOFK}, when $\Delta t_{k} \le 1/(2M_{B})$, we have
\begin{align}
     & \left|\rho_{k} - 1\right| =  \left|\frac{(f(x_{k}) - f(x_{k}+Pd_{k}))
     - (q_{k}(x_{k}) - q_{k}(x_{k}+Pd_{k}))}{q_{k}(x_{k}) - q_{k}(x_{k}+Pd_{k})}\right|
     \nonumber \\
     & \quad = \left|\frac{0.5(Pd_{k})^{T}\left(B_{k} - \nabla^2 f(\tilde{x}_{k})\right)(Pd_{k})}
     {q_{k}(x_{k}) - q_{k}(x_{k}+Pd_{k})}\right|
     \le \frac{0.5(M_{B} + M_{G} + L_{c} \|Pd_{k}\|) \|Pd_{k}\|^2}{|q_{k}(x_{k}) - q_{k}(x_{k}+Pd_{k})|}
     \nonumber \\
    & \quad \le \frac{(M_{B} + M_{G} + L_{c} \|Pd_{k}\|) \|Pd_{k}\|^2}
    {\left\|p_{g_{k}} \right\|
    \min\left\{\left\|Pd_{k}\right\|, \; \|p_{g_{k}}\|/(3\|B_{k})\|\right\}}
    \le \frac{(M_{B} + M_{G} + L_{c} \|Pd_{k}\|)\|Pd_{k}\|^2}
    {\epsilon_{p_g} \min\left\{\left\|Pd_{k}\right\|, \; \epsilon_{p_g}/(3M_{B})\right\}},
    \label{ESTRHOK}
\end{align}
where the last inequality is obtained from the assumption \eqref{PGKGEPN} of $p_{g_{k}}$.
We denote
\begin{align}
     M_{d} \triangleq \min\left\{\frac{\eta_{1}\epsilon_{p_{g}}}{M_{B}
    + M_{G}+L_{c}\epsilon_{p_g}/(3M_{B})}, \; \frac{\epsilon_{p_{g}}}{3M_{B}}\right\}.
    \label{UPBMPD}
\end{align}
Then, from equation \eqref{ESTRHOK}-\eqref{UPBMPD}, when $\|Pd_{k}\| \le M_{d}$, it is not
difficult to verify
\begin{align}
    \left|\rho_{k} - 1\right| \le \eta_{1}.     \label{RHOLETA1}
\end{align}

\vskip 2mm

From equations \eqref{PRDK}, \eqref{PROPD}, \eqref{PGGKBD} and \eqref{SEGLBD},
when $\Delta t_{k} \le 1/(2M_{B})$, we have
\begin{align}
    & \|Pd_{k}\|  \le \left\|d_{k}\right\| = \left\|\left(1/{\Delta t_{k}}I
    + B_{k} \right)^{-1}p_{g_{k}}\right\| \le {\|p_{g_{k}}\|}/(1/\Delta t_{k} - \|B_{k}\|) \nonumber \\
    & \quad \le {\|g_{k}\|}/(1/\Delta t_{k} - \|B_{k}\|)
    \le {M_g}/(1/\Delta t_{k} - M_{B}). \label{DKLECON}
\end{align}
We denote
\begin{align}
     \delta_{\Delta t} \triangleq \min\left\{ M_{d}/(M_g + M_{d}M_{B}), \;
     1/(2M_{B}) \right\}. \label{DELTAKUPBD}
\end{align}
Thus, from equations \eqref{DKLECON}-\eqref{DELTAKUPBD}, when $\Delta t_{k} \le
\delta_{\Delta t}$, we have $\|Pd_{k}\| \le M_{d}$.  That is to say, the condition of
inequality \eqref{RHOLETA1} holds.

\vskip 2mm

We assume that $K$ is the first index such that $\Delta t_{K} \le
\delta_{\Delta t}$ where $\delta_{\Delta t}$ is defined by equation \eqref{DELTAKUPBD}.
Then, from equations \eqref{RHOLETA1}-\eqref{DELTAKUPBD}, we know that
$|\rho_{K} - 1 | \le \eta_{1}$. According to the time-stepping adjustment
formula \eqref{ADTK1}, $x_{K} + Pd_{K}$ will be accepted and the time-stepping size
$\Delta t_{K+1}$ will be enlarged. Consequently, the time-stepping size $\Delta t_{k}$ holds
$\Delta t_{k}\ge \gamma_{2}\delta_{\Delta t}$ for all $k = 1, \, 2, \ldots$. \qed

\vskip 2mm

By using the results of Lemma \ref{LBSOAM} and Lemma \ref{DTBOUND}, we prove
the global convergence of Algorithm \ref{ALGPTCTR} for the linearly
constrained optimization problem \eqref{LEQOPT} as follows.

\vskip 2mm

\begin{theorem}
Assume that $f: \; \Re^{n} \rightarrow \Re$ is twice continuously differentiable
and $\nabla^{2} f(\cdot)$ is Lipschitz continuous. Moreover, we assume that the
level set $S_{f}$ defined by equation \eqref{LSCONFBD} and the quasi-Newton
matrices $B_{k}\, (k=1,\,2,\, \ldots)$ are bounded. The sequence $\{x_{k}\}$ is
generated by Algorithm \ref{ALGPTCTR}. Then, we have
\begin{align}
  \lim_{k \to \infty} \inf \|p_{g_{k}}\| = 0, \nonumber
\end{align}
where $p_{g_{k}} = P\nabla f(x_{k})$ and matrix $P$ is defined by equation \eqref{PROMAT}.
\end{theorem}
\proof We prove this result by contradiction as follows. Assume that the
conclusion is not true. Then, it exists a positive constant $\epsilon_{p_{g}}$
such that
\begin{align}
  \|p_{g_{k}}\| \ge \epsilon_{p_{g}} > 0, \;  k = 1, \, 2, \ldots.
  \label{PGKGECON}
\end{align}
According to Lemma  \ref{DTBOUND}, we know that it exists
an infinite subsequence $\{x_{k_{i}}\}$ such that trial steps $Pd_{k_i}$
are accepted, i.e., $\rho_{k_{i}} \ge \eta_{a}, \, i=1,\, 2,\ldots$. Otherwise,
all steps are rejected after a given iteration index, then the time-stepping
size will keep decreasing, which contradicts \eqref{DTGEPN}. Therefore,
from equation \eqref{MRHOK}, we have
\begin{align}
    f_{0} - \lim_{k \to \infty} f_{k} = \sum_{k = 0}^{\infty} (f_{k} - f_{k+1})
    \ge \eta_{a} \sum_{i = 0}^{\infty}
    \left(q_{k_{i}}(x_{k_i}) - q_{k_{i}}(x_{k_i}+Pd_{k_i})\right),
    \label{LIMSUMFK}
\end{align}
where $d_{k_{i}}$ is computed by equation \eqref{PRDK}.

\vskip 2mm

From the bounded assumption of $f(x)$ on the level set $S_{f}$ and
equation \eqref{LIMSUMFK}, we have
\begin{align}
    \lim_{k_{i} \to \infty}
    \left(q_{k_{i}}(x_{k_i}) - q_{k_{i}}(x_{k_i}+Pd_{k_i})\right) = 0.
    \label{LIMQK}
\end{align}
By substituting the estimate \eqref{PLBREDST} into equation \eqref{LIMQK}, we
obtain
\begin{align}
     \lim_{k_{i} \to \infty} \|p_{g_{k_i}}\|
     \min\left\{\left\|Pd_{k_i}\right\|, \;
     {\|p_{g_{k_i}}\|}/\left(3\|B_{k_i}\|\right)\right\} = 0.  \label{LIMPK}
\end{align}
According to the bounded assumptions of the level set $S_{f}$ and the
quasi-Newton matrices $B_{k}\, (k=1,\,2, \, \ldots)$, it exists two
positive constants $M_{p_{g}}$ and $M_{B}$ such that
\begin{align}
   \|p_{g_{k}}\| \le M_{p_{g}}, \; \|B_{k}\| \le M_{B}, \; k = 1, \, 2, \, \ldots.
   \label{PGGKBD2}
\end{align}
By substituting the bounded assumption \eqref{PGKGECON} of $p_{g_k} \,
(k = 1, \, 2, \, \ldots)$ and the bounded assumption \eqref{PGGKBD2} of matrices
$B_{k} \, (k = 1, \, 2, \, \ldots)$ into equation \eqref{LIMPK}, we have
\begin{align}
    \lim_{k_{i} \to \infty} \|Pd_{k_i}\| = 0.   \label{LIMPDK}
\end{align}

\vskip 2mm

From the Lipschitz continuous assumption \eqref{LIPSCHCON} of $\nabla^{2} f(\cdot)$, the
bounded assumption \eqref{PGGKBD2} of matrices $B_{k} \, (k = 1, \, 2, \, \ldots)$
and the bounded assumption of the level set $S_{f}$, we know that the result
\eqref{DTGEPN} of Lemma \ref{DTBOUND} is true. That is to say, it exists a positive
constant $\delta_{\Delta t}$ such that
\begin{align}
     \Delta t_{k} \ge \gamma_{2} \delta_{\Delta t} > 0, \;  k = 1,\, 2, \ldots.
 \label{DELTKGE}
\end{align}

\vskip 2mm

From equation \eqref{PRDK} and the property $P^{2} = P$ of the projection matrix $P$,
we obtain
\begin{align}
     Pd_{k_i} = \left(P \left({1}/{\Delta t_{k_i}}I
     + B_{k_i}\right)^{-1}P\right)p_{g_{k_i}}. \label{PDKEQ}
\end{align}
Thus, we have
\begin{align}
      & p_{g_{k_i}}^{T}Pd_{k_i}  = p_{g_{k_i}}^{T}\left(P \left({1}/{\Delta t_{k_i}}I
      + B_{k_i}\right)^{-1}P\right)p_{g_{k_i}} \nonumber \\
      & \quad = p_{g_{k_i}}^{T}\left({1}/{\Delta t_{k_i}}I
       + B_{k_i}\right)^{-1}p_{g_{k_i}} \ge \|p_{g_{k_i}}\|^2
       /(1/(\gamma_{2}\delta_{\Delta t}) + M_{B}). \label{PGKDK}
\end{align}
By applying the Cauchy-Schwartz inequality $|x^{T}y| \le \|x\|\|y\|$ to inequality
\eqref{PGKDK}, we have
\begin{align}
     \|p_{g_{k_i}}\|^2/(1/(\gamma_{2}\delta_{\Delta t}) + M_{B}) \le
     |p_{g_{k_i}}^{T}Pd_{k_i}| \le \|p_{g_{k_i}}\| \|Pd_{k_i}\|. \nonumber
\end{align}
That is to say, we have
\begin{align}
      \|p_{g_{k_i}}\| \le \left({1}/(\gamma_{2}\delta_{\Delta t}) + M_{B}\right)
      \|Pd_{k_i}\|. \label{PGKLEPDK}
\end{align}
By substituting the estimate \eqref{PGKLEPDK} into equation \eqref{LIMPDK},
we obtain
\begin{align}
     \lim_{{k_i} \to \infty} \|p_{g_{k_i}}\| = 0,  \label{PGKTOZ}
\end{align}
which contradicts the bounded assumption \eqref{PGKGECON} of
$p_{g_k} \, (k = 1, \, 2, \, \ldots)$.  \qed

\section{Numerical Experiments}

\vskip 2mm

In this section, some numerical experiments are executed to test the performance
of Algorithm \ref{ALGPTCTR} (the Ptctr method). The codes are performed by a Dell
G3 notebook with the Intel quad-core CPU and 8G memory. In subsection \ref{SubsecSans},
we compare Ptctr with the traditional optimization methods, i.e. the penalty
function method (PFM) \cite{Bertsekas2018} and SQP (the built-in subroutine
fmincon.m of the MATLAB2018a environment) \cite{FP1963,Goldfarb1970,MATLAB,NW1999},
and the traditional dynamical method, i.e. the backward differentiation formulas
(the built-in subroutine ode15s.m of the MATLAB2018a environment \cite{MATLAB,Shampine2002}),
for some large-scale linearly constrained-equality optimization problems which
are listed in Appendix A.

\vskip 2mm

According to the numerical results of Figures \ref{fig:CNMTIM}-\ref{fig:CNMTIM_LS}
and Tables \ref{TABCOM}- \ref{TABCOM_LS}, we find that Ptctr is superior to the other
three methods. Especially, the consumed time of Ptctr is much less than that of the
other three methods. In order to verify the performance of Ptctr further, we apply
it to a real-world optimization problem which arises from the visual-inertial
navigation localization when the unmanned aerial vehicle maintains the horizontal
flight, and compare it with SQP and the traditional dynamical method in subsection
\ref{SubSecPraProb}.

\vskip 2mm

For fairness, we use the pre-treatments of the constraints and the initial point
before we call the compared methods to solve the optimization problem
\eqref{LEQOPT} in the following test problems. The pre-processing methods are
stated \textbf{in subsection \ref{SUBSECDEG}}.

\vskip 2mm

\subsection{Statistical Analysis of Numerical Results} \label{SubsecSans}

\vskip 2mm

Here, the penalty factors $\sigma_{k} \, (k = 1, \, 2, \, \ldots)$ of PFM
\eqref{PENFX} are selected as $\sigma_{k+1} = 10 \sigma_{k} \; (k = 1, \, 2, \, \ldots)$
for the sequential unconstrained optimization subproblem. The subproblem
$\min_{x \in \Re^{n}} P_{\sigma_{k}}(x)$ is solved by the built-in subroutine
fminunc.m of the MATLAB2018a environment and its Hessian matrices are updated
by the BFGS method (pp. 194-198, \cite{NW1999}). We set the termination condition
of fminunc.m as
\begin{align}
       \|P_{\sigma_{k}}(x)\|_{\infty} \, \le 1.0 \times 10^{-8}, \;
        k = 1, \, 2, \, \ldots.      \nonumber
\end{align}
The initial point of the subproblem is the optimal solution of the previous
subproblem.

\vskip 2mm

We adopt ode15s.m to solve the ODE \eqref{ODGF} on the time interval
$[0, \; \infty)$ as the compared traditional dynamical method. Its relative
tolerance and absolute tolerance are both set by $1.0 \times 10^{-6}$.

\vskip 2mm

The termination conditions of the four compared methods are all set by
\begin{align}
    & \|\nabla_{x} L(x_{k}, \, \lambda_{k})\|_{\infty} \le 1.0 \times 10^{-6},
    \label{FOOPTTOL} \\
   & \|Ax_k - b \|_{\infty} \le 1.0 \times 10^{-6}, \;
    k = 1, \, 2, \, \ldots,   \label{FEATOL}
\end{align}
where the Lagrange function $L(x, \, \lambda)$ is defined by equation \eqref{LAGFUN}
and $\lambda$ is defined by equation \eqref{LAMBDA}.

\vskip 2mm

\begin{figure}[!htbp]
  \begin{minipage}[t]{0.48\linewidth}
      \centering
        \includegraphics[width= 1 \textwidth,height=0.25 \textheight]{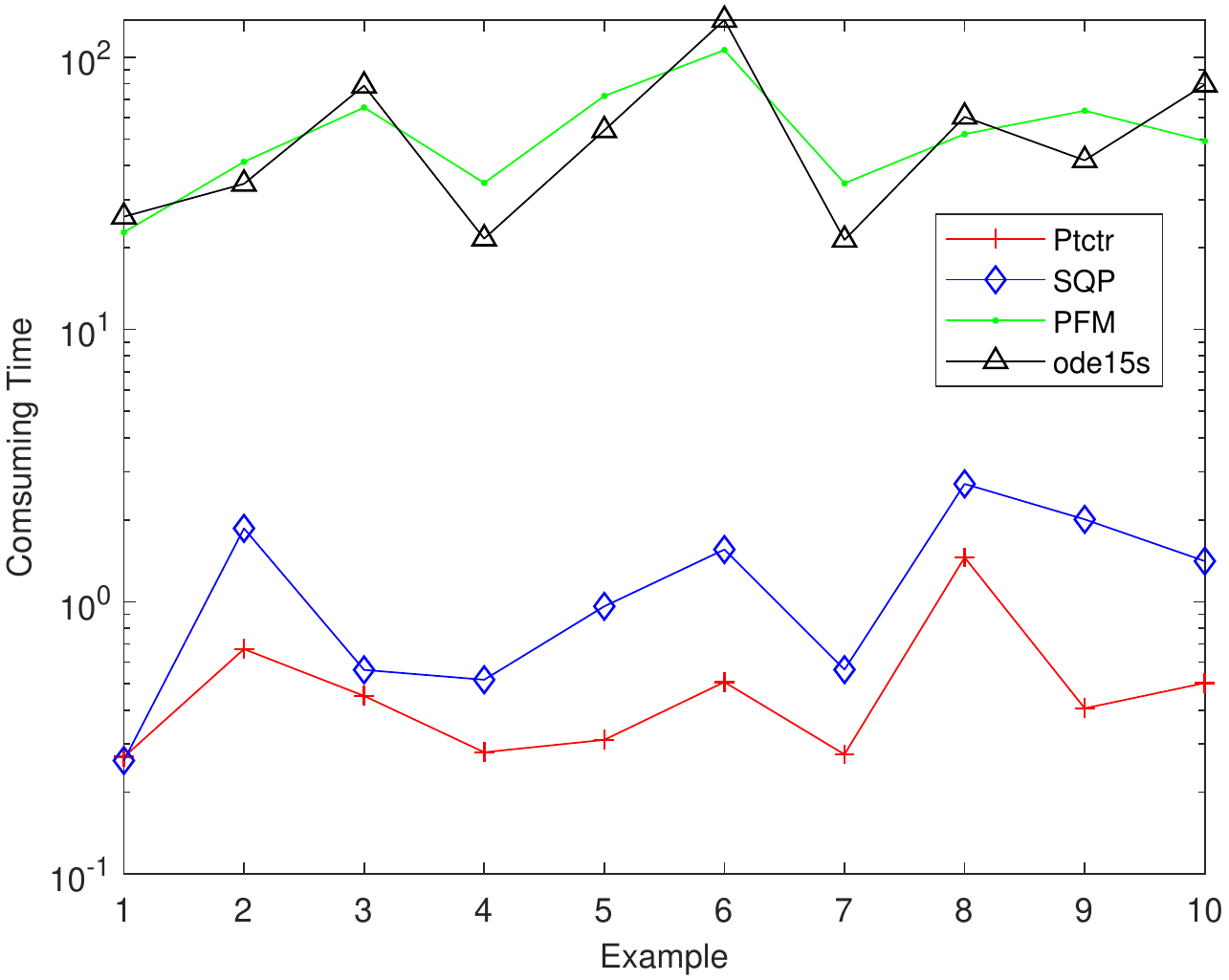}
        \caption{The consumed CPU time (s) of Ptctr, \, SQP, \, PFM, \, ode15s for test problems with $n \approx 1000$.}
        \label{fig:CNMTIM}
  \end{minipage}
  \hfill
  \begin{minipage}[t]{0.48\linewidth}
      \centering
            \includegraphics[width= 1 \textwidth,height=0.25\textheight]{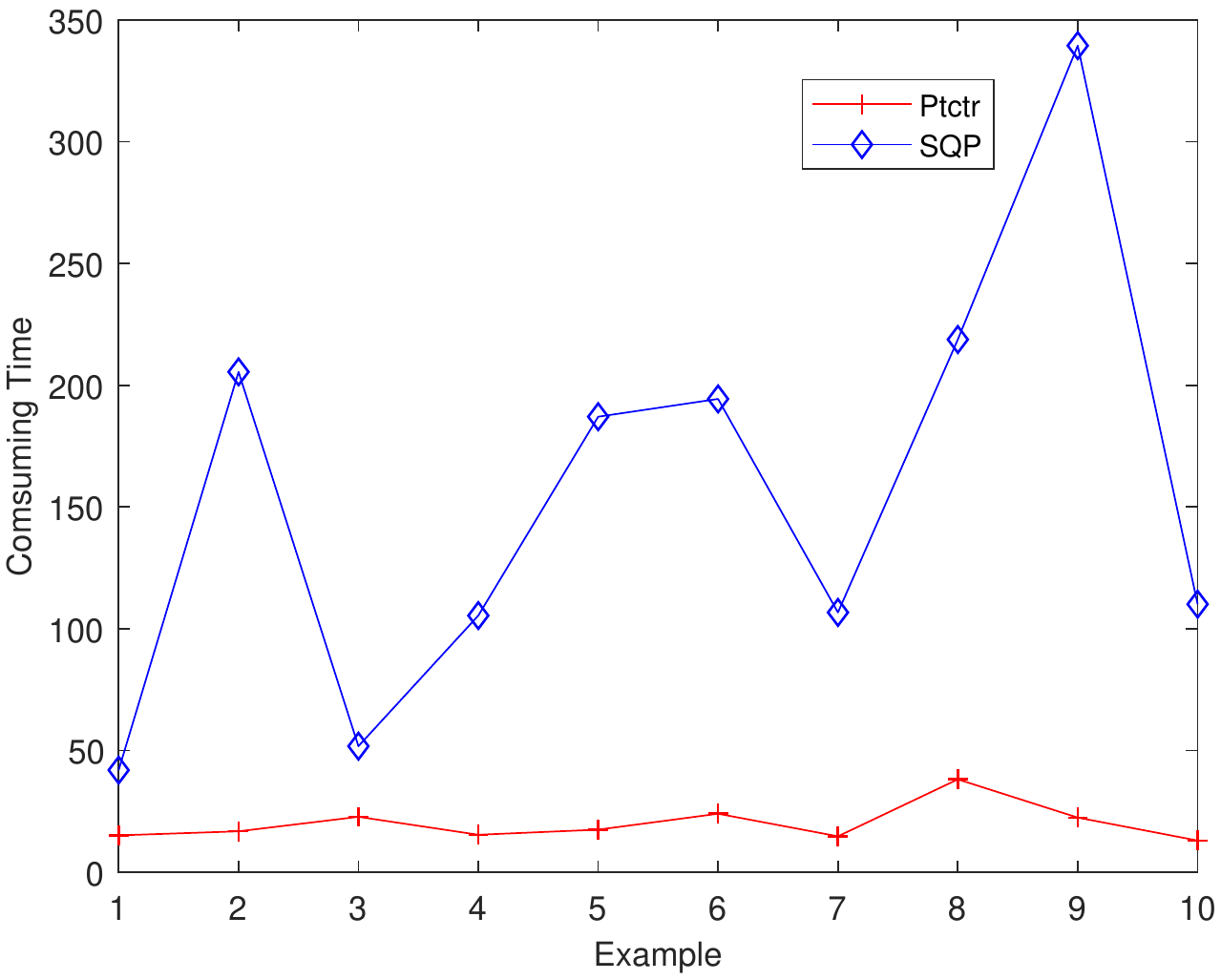}
            \caption{The consumed CPU time (s) of Ptctr, \, SQP for test problems with $n \approx 5000$.}
            \label{fig:CNMTIM_LS}
  \end{minipage}%
\end{figure}

Ten large-scale test problems are listed in Appendix A. \textbf{Firstly, we
test one set of data for those ten problems with $n \approx 1000$. The numerical
results are put in Table \ref{TABCOM} and Figure \ref{fig:CNMTIM}. From Table
\ref{TABCOM}, we find that Ptctr and SQP can correctly solve those ten test
problems with $n \approx 1000$. However, PFM only approaches the correct solutions
of eight test problems since it can not attain the KKT termination condition
\eqref{FOOPTTOL} of those test problems. From Table \ref{TABCOM} and Figure
\ref{fig:CNMTIM}, we find that the consumed time of Ptctr and SQP are far less
than those of the other two methods, respectively.}

\vskip 2mm

\textbf{In order to discriminate the performances of Ptctr and SQP further, we
test another set of data for these ten problems with $n \approx 5000$. The
numerical results are put in Table \ref{TABCOM_LS} and Figure \ref{fig:CNMTIM_LS}.
From Table \ref{TABCOM_LS} and Figure \ref{fig:CNMTIM_LS}, we find that Ptctr and
SQP can correctly solve those ten test problems with $n \approx 5000$ and the consumed
time of Ptctr is significantly less than that of SQP for every test problem (the consumed 
time of Ptctr is about one fifth of that of SQP for the non-quadratic programming problem).}

\vskip 2mm

\textbf{From those test sets of data, we find that Ptctr works significantly better
than the other three methods, respectively. One of reasons is that Ptctr only solves
a linear system of equations with an $n \times n$ symmetric definite coefficient matrix
at every iteration and it requires about $n^{3}/3$ flops since we use the Cholesky
factorization to solve it (p. 169, \cite{GV2013}). However, SQP needs to solve a linear
system of equations with dimension $(m+n)$ when it solves a quadratic programming subproblem
at every iteration (pp. 531-532, \cite{NW1999}) and it requires about $2(m+n)^{3}/3$ flops
(p. 116, \cite{GV2013}). The other two methods need to solve a nonlinear system of equations
(ode15s), or an unconstrained optimization subproblem (PFM) at every iteration. }

\vskip 2mm

\begin{table}[!http]
  \newcommand{\tabincell}[2]{\begin{tabular}{@{}#1@{}}#2\end{tabular}}
  \scriptsize
  \centering
  \caption{Numerical results of test problems with $n \approx 1000$.} \label{TABCOM}
  \begin{tabular}{|c|c|c|c|c|c|c|c|c|}
  \hline
  \multirow{2}{*}{Problems } & \multicolumn{2}{c|}{Ptctr} & \multicolumn{2}{c|}{SQP} & \multicolumn{2}{c|}{PFM} & \multicolumn{2}{c|}{ode15s} \\ \cline{2-9}
                        & \tabincell{c}{steps \\(time)}     & $f(x^\star)$    & \tabincell{c}{steps \\(time)}  & $f(x^\star)$   & \tabincell{c}{steps \\(time)}& $f(x^\star)$ & \tabincell{c}{steps \\(time)} & $f(x^\star)$

                        \\ \hline
  \tabincell{c}{Exam. 1 \\ (n = 1000, \\ m = n/2)}   &\tabincell{c}{11 \\ (0.27)} &7.27E+03   &\tabincell{c}{2 \\ (0.26)} &7.27E+03   &\tabincell{c}{18 \\ (22.75)} & \tabincell{c}{7.27E+03 \\ (close)}  &\tabincell{c}{81 \\ (25.96)} &7.27E+03   \\ \hline
  \tabincell{c}{Exam. 2 \\ (n = 1000,\\ m = n/3)}   &\tabincell{c}{18 \\ (0.67)} &1.29E+03   &\tabincell{c}{17 \\ (1.86)} &1.29E+03   &\tabincell{c}{14 \\ (41.28)} & \tabincell{c}{1.29E+03 \\ (close)} &\tabincell{c}{123 \\ (34.27)} & 1.29E+03    \\ \hline
  \tabincell{c}{Exam. 3 \\ (n = 1200, \\m = 2n/3)}   &\tabincell{c}{12 \\ (0.45)} &714.67     &\tabincell{c}{2 \\ (0.56)} &714.67     &\tabincell{c}{15\\ (65.36)} &714.67           &\tabincell{c}{66 \\ (78.57)} &714.67  \\\hline
  \tabincell{c}{Exam. 4 \\ (n = 1000, \\m = n/2)} 	&\tabincell{c}{11 \\ (0.28)}	&97.96	    &\tabincell{c}{6 \\ (0.52)} &97.96      &\tabincell{c}{15 \\ (34.53)} &\tabincell{c}{97.96 \\ (close)}   &\tabincell{c}{62 \\ (21.56)}	&97.96  \\ \hline
  \tabincell{c}{Exam. 5 \\ (n = 1000, \\m = n/2)} 	&\tabincell{c}{14 \\ (0.31)}&82.43	    &\tabincell{c}{11 \\ (0.96)}	&82.43      &\tabincell{c}{15 \\ (72.05)} &\tabincell{c}{82.53 \\ (close)}    &\tabincell{c}{173\\ (53.88)}	&82.43  \\ \hline
  \tabincell{c}{Exam. 6 \\ (n = 1200, \\m = 2n/3)} 	&\tabincell{c}{13 \\ (0.51)} &514.48	    &\tabincell{c}{9 \\ (1.56)}	&514.48     &\tabincell{c}{15 \\ (106.19)} &\tabincell{c}{514.48  \\ (close)}    &\tabincell{c}{66 \\ (136.97)} &514.48   \\ \hline
  \tabincell{c}{Exam. 7 \\ (n = 1000, \\ m = n/2)} 	&\tabincell{c}{10 \\ (0.28)}  &1.19E+04	&\tabincell{c}{6 \\ (0.56)}	& 1.19E+04	&\tabincell{c}{16 \\ (34.43)} &\tabincell{c}{1.19E+04 \\ (close)}   &\tabincell{c}{77 \\ (21.34)} &1.19E+04\\ \hline
  \tabincell{c}{Exam. 8 \\ (n = 1200, \\m = n/3)} 	&\tabincell{c}{38 \\ (1.46)} &196.24	    &\tabincell{c}{26\\(2.71)}	&196.24	    &\tabincell{c}{16 \\ (52.20)}&\tabincell{c}{197.83  \\ (close)}    &\tabincell{c}{147 \\ (60.41)} &196.24 \\ \hline
  \tabincell{c}{Exam. 9 \\ (n = 1000, \\ m = n/2)} 	&\tabincell{c}{13 \\ (0.41)} &4.42E+04	&\tabincell{c}{29 \\ (2.01)}	&4.42E+04	&\tabincell{c}{17 \\ (63.61)} &\tabincell{c}{4.42E+04  \\ (close)}   &\tabincell{c}{118 \\ (41.81)} &4.42E+04  \\ \hline
  \tabincell{c}{Exam. 10 \\ (n = 1200, \\ m = n/3)} 	&\tabincell{c}{16 \\ (0.50)} &0.50      &\tabincell{c}{14 \\ (1.41)}	&0.50 &\tabincell{c}{11 \\ (49.17)}	&0.50   &\tabincell{c}{179 \\ (79.25)} &0.50\\ \hline
\end{tabular}
\end{table}

\vskip 2mm

\begin{table}[!http]
  \newcommand{\tabincell}[2]{\begin{tabular}{@{}#1@{}}#2\end{tabular}}
  \scriptsize
  \centering
  \caption{Numerical results of test problems with $n \approx 5000$.} \label{TABCOM_LS}
  \begin{tabular}{|c|c|c|c|c|}
  \hline
  \multirow{2}{*}{Problems } & \multicolumn{2}{c|}{Ptctr} & \multicolumn{2}{c|}{SQP} \\ \cline{2-5}
                        & \tabincell{c}{steps \\(time)}     & $f(x^\star)$    & \tabincell{c}{steps \\(time)}  & $f(x^\star)$

                        \\ \hline
  \tabincell{c}{Exam. 1 \\ (n = 5000, \\ m = n/2)}   &\tabincell{c}{11 \\ (15.17)} &3.636364E+04   &\tabincell{c}{2 \\ (41.96)} &3.636364E+04      \\ \hline
  \tabincell{c}{Exam. 2 \\ (n = 5000,\\ m = n/3)}   &\tabincell{c}{16 \\ (16.84)} &5.179806E+03   &\tabincell{c}{6 \\ (205.50)} &5.179806E+03  \\ \hline
  \tabincell{c}{Exam. 3 \\ (n = 4800, \\m = 2n/3)}   &\tabincell{c}{12 \\ (22.80)} &2.858667E+03     &\tabincell{c}{2 \\ (51.79)} &2.858667E+03     \\\hline
  \tabincell{c}{Exam. 4 \\ (n = 5000, \\m = n/2)} 	&\tabincell{c}{11 \\ (15.42)}	&4.937947E+02	    &\tabincell{c}{6 \\ (105.43)} &4.937947E+02      \\ \hline
  \tabincell{c}{Exam. 5 \\ (n = 5000, \\m = n/2)} 	&\tabincell{c}{14 \\ (17.53)} &4.321521E+02	    &\tabincell{c}{11 \\ (187.08)} & 4.321521E+02      \\ \hline
  \tabincell{c}{Exam. 6 \\ (n = 4800, \\m = 2n/3)} 	&\tabincell{c}{13 \\ (24.06)} & 2.057906E+03	    &\tabincell{c}{11 \\ (194.38)}	&2.057906E+03     \\ \hline
  \tabincell{c}{Exam. 7 \\ (n = 5000, \\m = n/2)} 	&\tabincell{c}{10 \\ (14.79)}  &5.944739E+04	&\tabincell{c}{6 \\ (106.72)}	& 5.944739E+04 	\\ \hline
  \tabincell{c}{Exam. 8 \\ (n = 4800, \\m = n/3)} 	&\tabincell{c}{38 \\ (38.17)} &7.768754E+02	    &\tabincell{c}{28\\(218.80)}	& 7.768754E+02	    \\ \hline
  \tabincell{c}{Exam. 9 \\ (n = 5000, \\m = n/2)} 	&\tabincell{c}{12 \\ (22.47)} &2.211073E+05	&\tabincell{c}{22 \\ (339.47)} &2.211073E+05	 \\ \hline
  \tabincell{c}{Exam. 10 \\ (n = 4800, \\m = n/3)} 	&\tabincell{c}{16 \\ (13.06)} &2.002622E+00      &\tabincell{c}{14 \\ (110.11)}	&2.002622E+00 \\ \hline
\end{tabular}
\end{table}

\vskip 2mm

\vskip 2mm

\subsection{The Visual-Inertial Navigation Localization Problems}
\label{SubSecPraProb}

\vskip 2mm

In this subsection, in order to verify the performance of Ptctr further,
we apply it to a real-world problem which arises from the visual-inertial
navigation localization problem when the unmanned aerial vehicle maintains
the horizontal flight, and compare it with SQP and ode15s.

\vskip 2mm

When an unmanned aerial vehicle flies with a speed of about 200 m/s at a
kilometer altitude for an hour, the positioning error of the pure inertial
navigation system is about ten \textbf{kilometers}, which cannot meet the requirement
of the positioning accuracy, i.e. less than one kilometer. Therefore, we
consider the visual-inertial navigation localization method for this problem
\cite{EBM2018,LLS2020}.

\vskip 2mm

\begin{figure}[htbp]
    \centering
    \includegraphics[width= 0.6 \textwidth,height=0.3 \textheight]{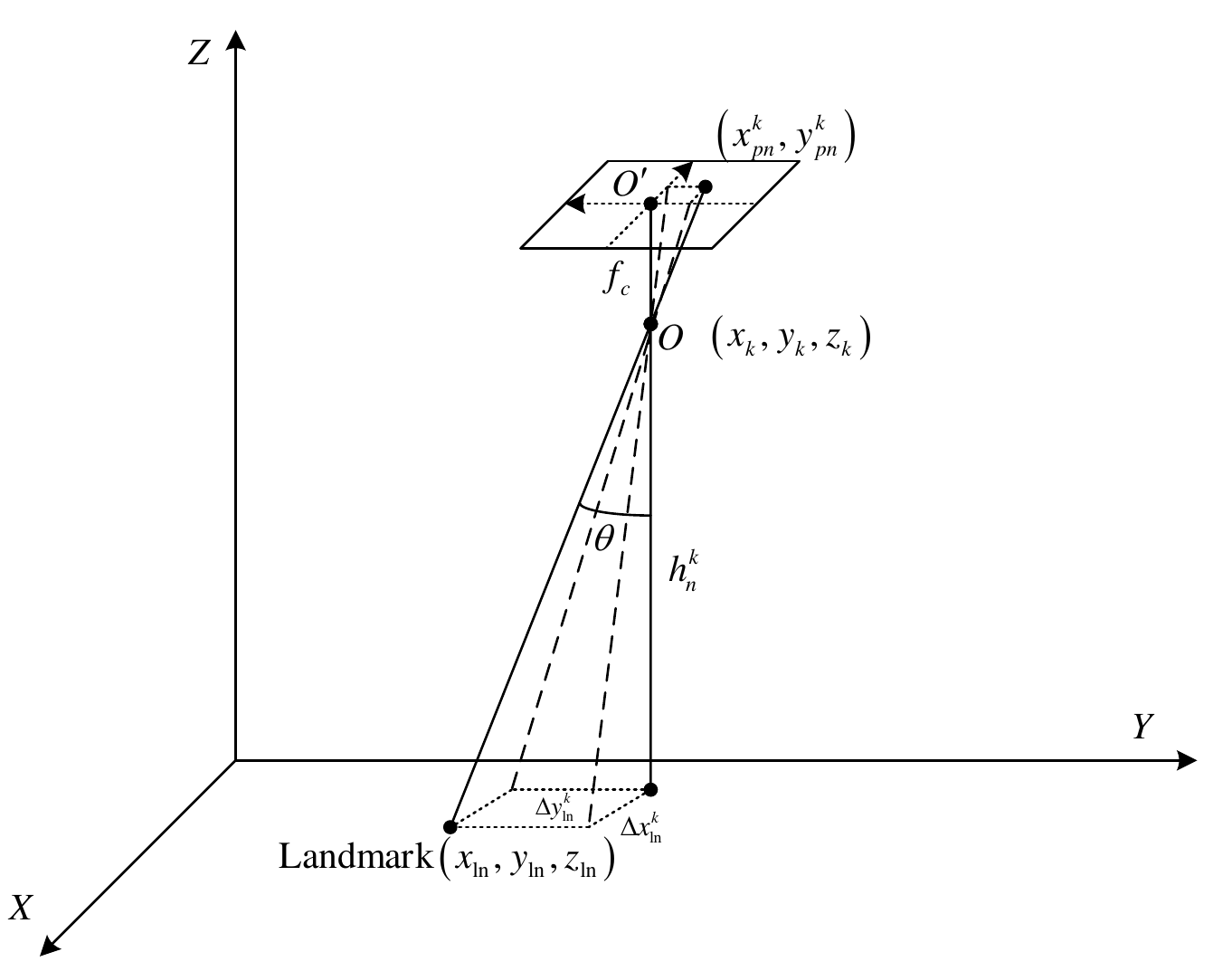}
    \caption{The pinhole camera model.}
    \label{Fig:PM}
\end{figure}

For the visual part, we convert the world coordinate system of landmarks to
the camera coordinate system by the pinhole camera model \cite{HZ2003}. The
pinhole model is illustrated in Figure \ref{Fig:PM}. There are two coordinate
systems, which are depicted specifically there. $O$ represents the optical
centre of the camera lens and $f_c$ refers to the focal length of the camera.
The $X^{'} $-axis and $Y{'}$ in the camera coordinate system are parallel
to the $X$-axis and the $Y$-axis in the world coordinate system, respectively.
The position of the  $k$-th camera in the world coordinate system is denoted
as $(x_k,\; y_k, \,z_k)$. $(x_{ln},\; y_{ln},\; z_{ln})$ denotes the position
of the $n$-th landmark in the world coordinate system. The vertical distance
between the $n$-th landmark and the camera position of the $k$-th frame is
denoted as $h_{n}^{k} = z_{k} - z_{ln}$. $\Delta{x_{n}^k},\,\Delta{y_{n}^k}$
respectively represent the $X$-axis and $Y$-axis coordinate differences between
the $n$-th landmark and the $k$-th camera in the world coordinate system.
$(x_{pn}^{k},\; y_{pn}^{k})$ represents the coordinate in the camera coordinate
system of the $k$-th frame, which is projected from the $n$-th landmark in the
world coordinate system. $\theta$ denotes the line-of-sight
angle of the landmark relative to the optical center of the camera.

\vskip 2mm

\vskip 2mm

Thus, by combining the visual information provided by the camera and the distance
information $dist_{hor}$ obtained from Inertial Measurement Units (IMU) between
the $k$-th frame camera and the $(k+1)$-th frame camera, we obtain the position
of $(x_{k+1}, \; y_{k+1})$ of the $(k+1)$-th frame camera in the world coordinate
system via solving the following $(k+1)$-th optimization problem:
\begin{align}
    \min_{(x_{k+1}, \; y_{k+1})} \quad &
    \left(\sqrt{(x_{k+1}-x_k)^2 + (y_{k+1} - y_{k})^{2}} - dist_{hor} \right)^2
    \label{DISFUN} \\
    \text{subject to} \; & x_{ln} + \frac{x_{pn}^{k}}{f_c}h_{n}^{k}  = x_{k}, \nonumber \\
        & y_{ln} + \frac{y_{pn}^{k}}{f_c}h_{n}^{k}  = y_{k}, \label{EQN:LIN} \\
        & x_{k+1} - x_{ln} - \frac{x_{pn}^{k+1}}{f_c}h_{n}^{k}
         = \frac{\Delta{h_{k}^{k+1}}}{f_c}x_{pn}^{k+1}, \nonumber \\
        & y_{k+1} - y_{ln} - \frac{y_{pn}^{k+1}}{f_c}h_{n}^{k}
         = \frac{\Delta{h_{k}^{k+1}}}{f_c}y_{pn}^{k+1}, \nonumber
\end{align}
where $(x_{k+1}, \, y_{k+1})$, $(x_{k}, \, y_{k})$ respectively represent the
$x$-$y$ axis coordinates of the camera at the $(k+1)$-th frame and the $k$-th frame in
the world coordinate system, $(x_{ln},\, y_{ln})$ represents the coordinate of the
$n$-th landmark in the world coordinate system, and $(x_{pn}^{k},\, y_{pn}^{k})$,
$(x_{pn}^{k+1},\, y_{pn}^{k+1})$ respectively represent the coordinates in the camera
coordinate system at the $k$-th frame and the $(k+1)$-th frame, which are projected
from the $n$-th landmark in the world coordinate system. In equation \eqref{EQN:LIN},
$\Delta{h_{k}^{k+1}}$ is the altitude difference obtained by an altimeter between
the $k$-th camera position and the $(k+1)$-th camera position, and $dist_{hor}$
represents the horizontal Euclidean distance provided by the IMU between the $k$-th
camera position and the $(k+1)$-th camera position.

\vskip 2mm

For the optimization problem \eqref{DISFUN}-\eqref{EQN:LIN}, $(x_{k+1}, \; y_{k+1})$,
$(x_{ln},\; y_{ln})$ and $\Delta{h_{k}^{k+1}}$ are unknown. Obviously, this is
an under-determined system if we only use one landmark and we can not uniquely
determine the position $(x_{k+1},\; y_{k+1})$ of the camera at the $(k+1)$-th
frame. Therefore, we use five landmarks to determine the position of the camera
at the $(k+1)$-th frame, and rearrange the constraint system \eqref{EQN:LIN} as
follows:
\begin{align}
      A_{k+1}w = b_{k+1}, \; A_{k+1} = \begin{bmatrix}
                     B &C_{1} &O &O &O &O \\
                     B &O &C_{2} &O &O &O\\
                     B &O&O &C_{3} &O &O \\
                     B &O &O &O &C_{4} &O \\
                     B &O &O &O &O &C_{5} \\
                   \end{bmatrix},
     \label{EQN:5LD}
\end{align}
where $O \in \Re^{4 \times 3}$ represents a zero matrix, and matrix $B$, matrices
$C_{n} \, (n = 1, \, 2, \, \ldots, \, 5)$, the constant vector $b_{k+1}$, and the
variable $w$  are defined as follows:
\begin{align}
      & \hskip 2cm B = \begin{bmatrix}
               0 &0 \\
               0 &0 \\
               1 &0 \\
               0 &1 \\
           \end{bmatrix},  \;
     C_{n} = \begin{bmatrix}
                  1 & 0 & \frac{x_{pn}^{k}}{f_c}  \\
                  0 &1 & \frac{y_{pn}^{k}}{f_c} \\
                  -1 &0 & -\frac{x_{pn}^{k+1}}{f_c} \\
                  0 &-1 & -\frac{y_{pn}^{k+1}}{f_c} \\
              \end{bmatrix}, \;  n = 1, \, 2, \, \ldots, \, 5,  \nonumber \\
     & b_{k+1} =\left[x_{k}, \,  y_{k}, \, \frac{\Delta{h_{k}^{k+1}}}{f_c}x_{p1}^{k+1},
      \, \frac{\Delta{h_{k}^{k+1}}}{f_c}y_{p1}^{k+1}, \, \ldots, \, x_{k}, \,
      y_{k}, \, \frac{\Delta{h_{k}^{k+1}}}{f_c}x_{p5}^{k+1}, \,
      \frac{\Delta{h_{k}^{k+1}}}{f_c}y_{p5}^{k+1} \right]^{T}, \nonumber \\
     & w = [x_{k+1}, \, y_{k+1}, \, x_{l1}, \,  y_{l1}, \, h_{1}^{k}, \, \ldots,
      \, x_{l5}, \, y_{l5}, \, h_{5}^{k}]^{T}. \nonumber
\end{align}

 \vskip 2mm

After simple calculations, we know that the rank of the linear system \eqref{EQN:5LD}
is deficient when  $\Delta h_{k}^{k+1} = 0$, where $\Delta h_{k}^{k+1}$ represents
the altitude difference between the $k$-th frame and the $(k+1)$-th frame of the
camera.

\vskip 2mm

Due to the measurement error for a real-world problem, the linear system \eqref{EQN:5LD}
has the following general form:
\begin{align}
   A_{k+1}(\varepsilon) w = b_{k+1}(\varepsilon), \label{RANDCON}
\end{align}
where $A_{k+1} (\varepsilon)$ includes the visual measurement error
$\varepsilon_{\theta}$ and the altimeter measurement error $\varepsilon_{h}$,
$b_{k+1} (\varepsilon)$ is the constant vector with noisy data.
Thus, the visual-inertial navigation localization problem \eqref{DISFUN}-\eqref{EQN:5LD}
is extended by the following linearly equality-constrained optimization problem
with noisy data:
\begin{align}
     \min_{w \in \Re^{17}} \quad & f_{k} (w) =
     \left(\sqrt{(x_{k+1}-x_k)^2 + (y_{k+1} - y_{k})^{2}}-dist_{hor,\varepsilon}\right)^2
     \label{ERRDISFUN} \\
     \text{subject to} \; & A_{k+1}(\varepsilon)s = b_{k+1}(\varepsilon),  \label{ERRCON}
\end{align}
where matrix $A_{k+1}(\varepsilon)$ and vector $b_{k+1}(\varepsilon)$ are defined by
equation \eqref{EQN:5LD}.

\vskip 2mm

After establishing the mathematical model \eqref{ERRDISFUN}-\eqref{ERRCON} of
the visual-inertial navigation localization problem, we simulate the trajectories
of the unmanned aerial vehicle maintaining the horizontal flight. By using the
optimization method to solve the linearly equality-constrained optimization problem
\eqref{ERRDISFUN}-\eqref{ERRCON} at every sampling time, we obtain the error propagation
of the unmanned aerial vehicle flying for one hour.

\vskip 2mm

According to the given condition from the industry, we assume that the unmanned
aerial vehicle flies horizontally for an hour at an altitude of 1200 meters with
speed 235 meters per second and the sampling period is half a second. The camera
focal length $f_{c}$ is set by $f_{c} = 24\times 10^{-3}$ meter. We simulate three
simple trajectories with or without measurement errors. The first trajectory is
represented as follows:
\begin{align}
     \text{trj 1:} \;  (x_{k}, \, y_{k}, \, z_{k}) =
     (0,  \, 117.5 \, k,\, 1200),
     \; k = 1, \, 2, \, \ldots, 7200. \label{TRA1}
\end{align}
The second trajectory is represented as follows:
\begin{align}
  &\text{trj 2:} \; (x_{k}, \, y_{k}, \, z_{k}) =
  \left(\frac{1}{2}d_{k}, \, \frac{\sqrt{3}}{2}d_{k}, \,  1200 \right), \;
  d_{k} = 117.5 \, k, \; k = 1, \, 2, \, \ldots, \, 1800,
  \nonumber \\
  &\text{and} \; (x_{k}, \, y_{k}, \, z_{k}) =
  \left(\frac{1}{2}d_{1800}, \, \frac{\sqrt{3}}{2}d_{1800} + d_{k}, \, 1200 \right), \;
  d_{k} = 117.5 \,  k, \nonumber \\
  & k = 1801, \, 1802, \, \ldots, \, 7200. \label{TRA2}
\end{align}
The third trajectory is represented as follows:
\begin{align}
   \text{trj 3:} \; (x_{k}, \, y_{k}, \, z_{k}) =
  \left(\frac{1}{2}d_{k}, \, \frac{\sqrt{3}}{2}d_{k}, \,  1200 \right), \;
  d_{k} = 117.5 \, k, \; k = 1, \, 2, \, \ldots, \, 7200.
  \label{TRA3}
\end{align}

\vskip 2mm

Thus, based on the known trajectory, we set the coordinate
$(x_{ln}, \, y_{ln}, \, z_{ln})_{k}$ of the $n$-th landmark observed by the $k$-th
camera in the world coordinate system as follows:
\begin{align}
   (x_{ln}, \, y_{ln}, \, z_{ln})_{k} = \left(x_{k} + \frac{58.75n}{N}, \,
   y_{k} + \frac{58.75n}{N}, \, \frac{40n}{N}\right),
    \; n = 1, \, 2, \, \ldots, \,  N.  \label{LANDMARK}
\end{align}
According to the principle of pinhole imaging, from equations \eqref{EQN:LIN}
and \eqref{LANDMARK}, we generate the image coordinate
$\left(x_{pn}^{k}, \, y_{pn}^{k}\right)$ of the $n$-th landmark
in the $k$-th camera coordinate system as follows:
\begin{align}
    x_{pn}^{k} = \frac{x_{k} - x_{ln}}{h_{n}^{k}} f_{c},
     \; y_{pn}^{k} = \frac{y_{k} - y_{ln}}{h_{n}^{k}} f_{c}, \; n = 1, \, 2,
     \ldots, \,  N,
     \label{PIXLANDMK}
\end{align}
where the altitude difference $h_{n}^{k}$ between the $k$-th camera and the
$n$-th landmark in the world coordinate system is computed by
$h_{n}^{k} = z_{k} - z_{ln}$. The altitude difference $\Delta{h_{k}^{k+1}}$
between the $k$-th camera position and the $(k+1)$-th camera position is
computed by $\Delta{h_{k}^{k+1}} = z_{k+1} - z_{k}$.

\vskip 2mm

\begin{figure}[!htbp]
   \begin{minipage}[t]{0.48\linewidth}
      \centering
      \subfigure[trajectory 1]{
        \includegraphics[width=1\textwidth,height=0.25\textheight]{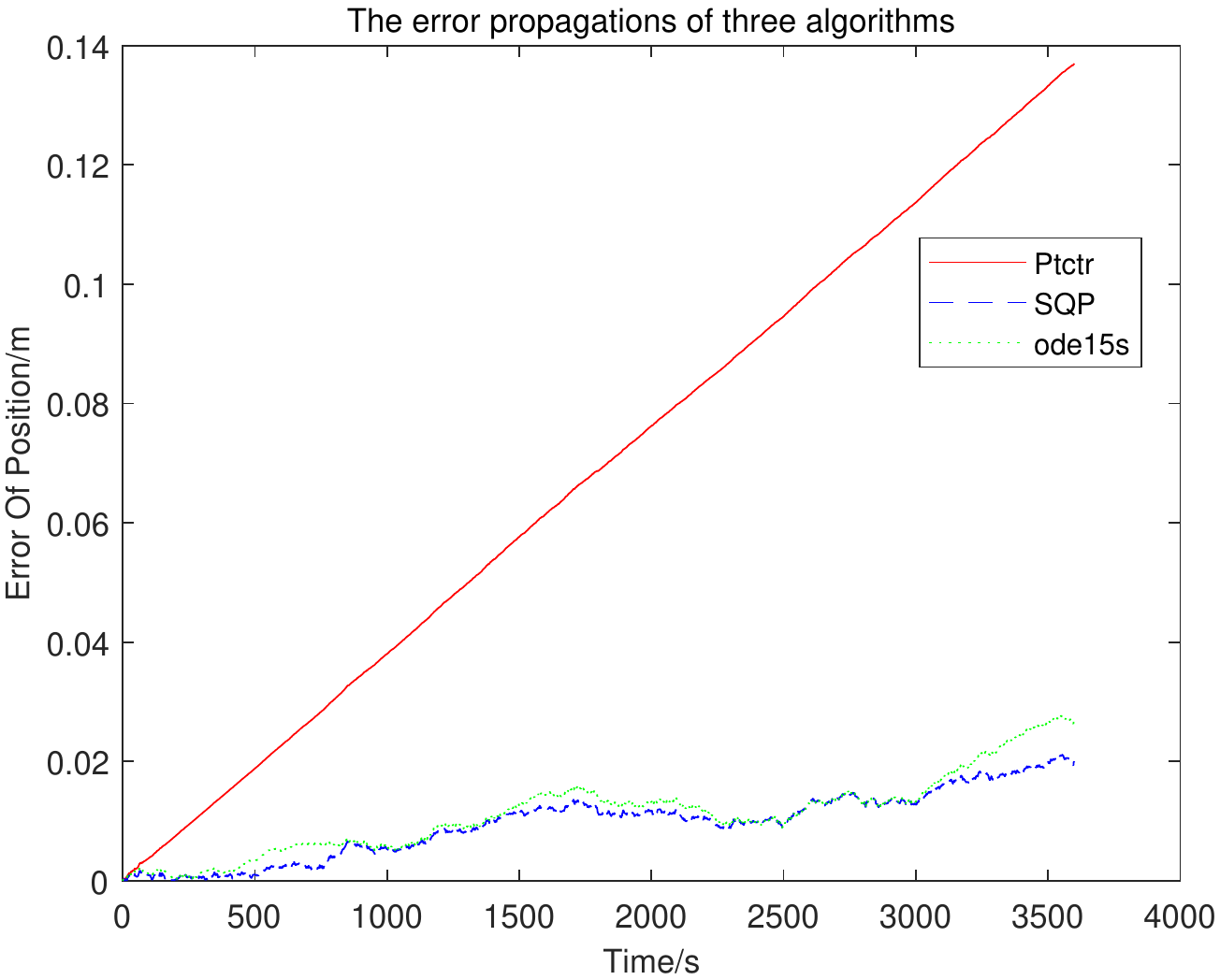}
      }
      \subfigure[trajectory 2]{
        \includegraphics[width=1\textwidth,height=0.25\textheight]{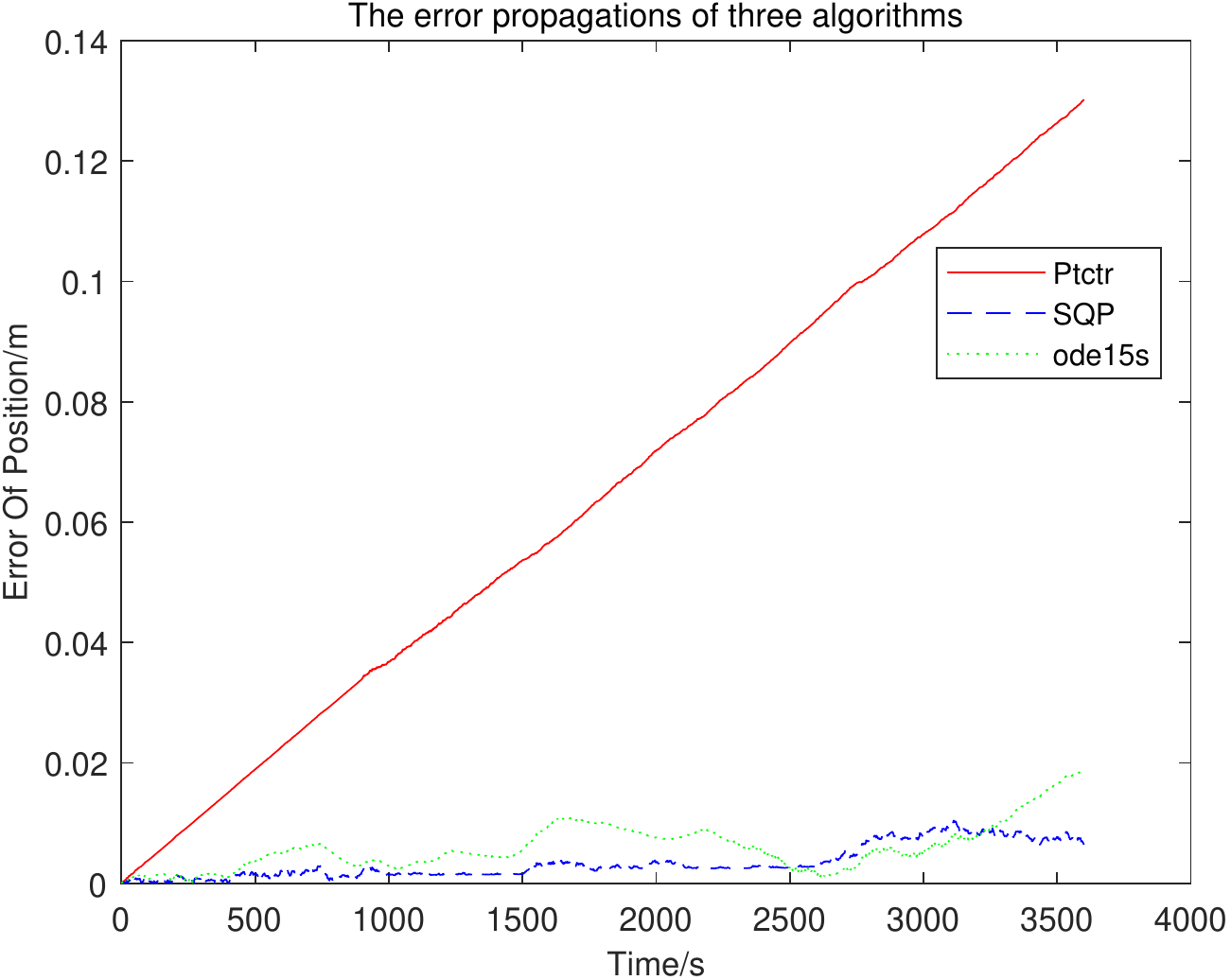}
      }
      \subfigure[trajectory 3]{
        \includegraphics[width=1\textwidth,height=0.25\textheight]{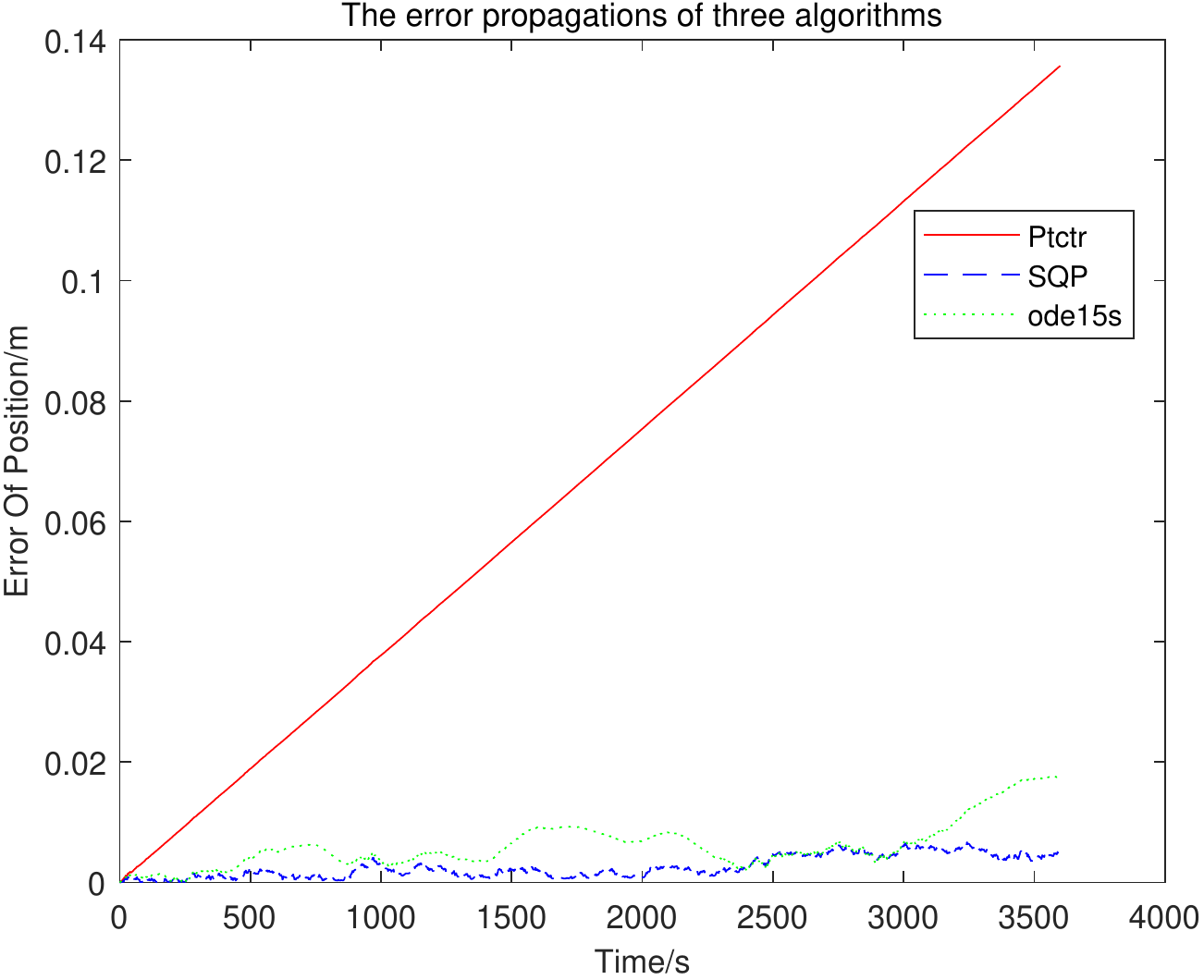}
      }
      \caption{The visual-inertial positioning \\
      problems without random errors.} \label{FIGNOERR}
   \end{minipage}%
   \hfill
   \begin{minipage}[t]{0.48\linewidth}
      \centering
      \subfigure[trajectory 1]{
        \includegraphics[width=1\textwidth,height=0.25\textheight]{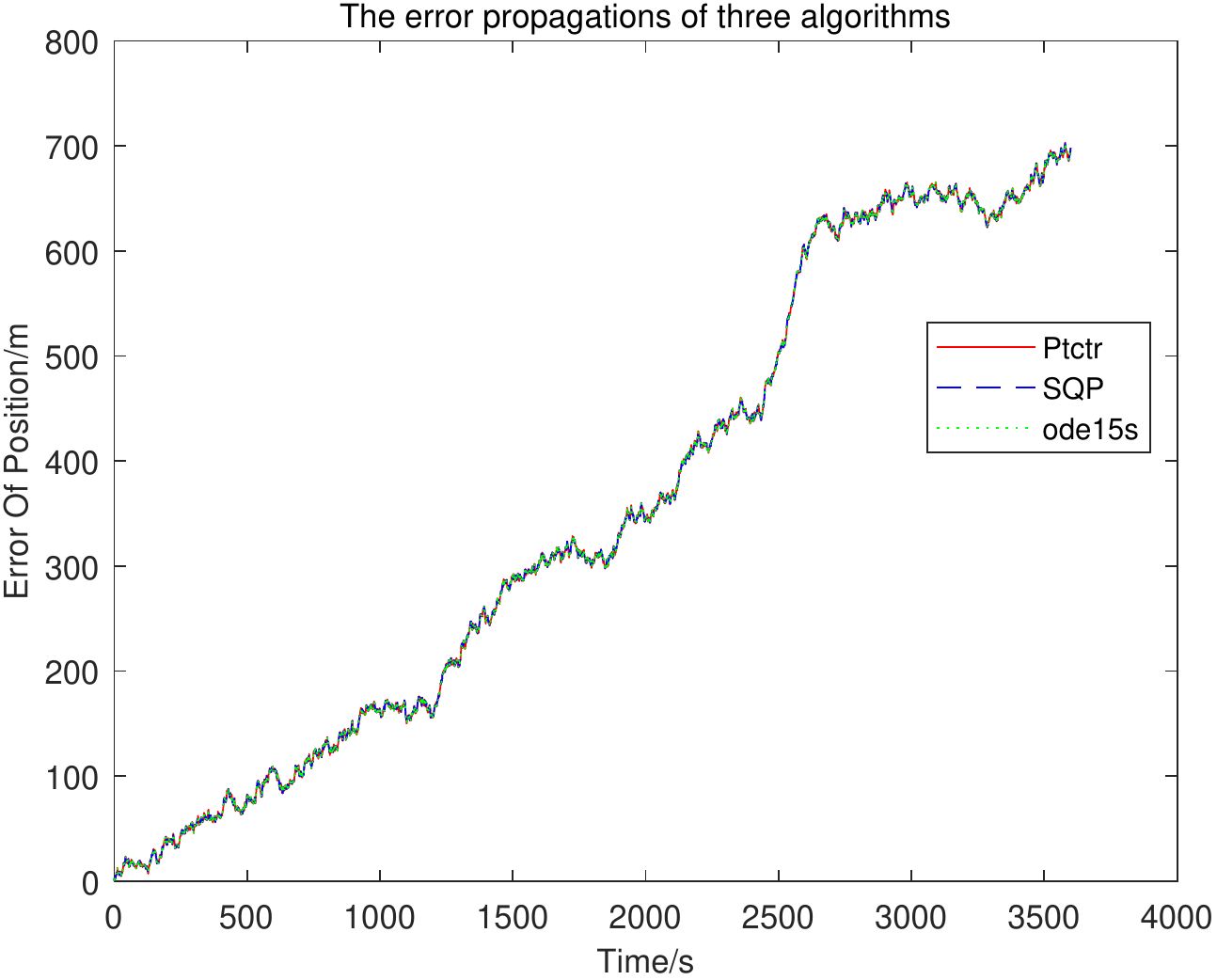}
      }
      \subfigure[trajectory 2]{
        \includegraphics[width=1\textwidth,height=0.25\textheight]{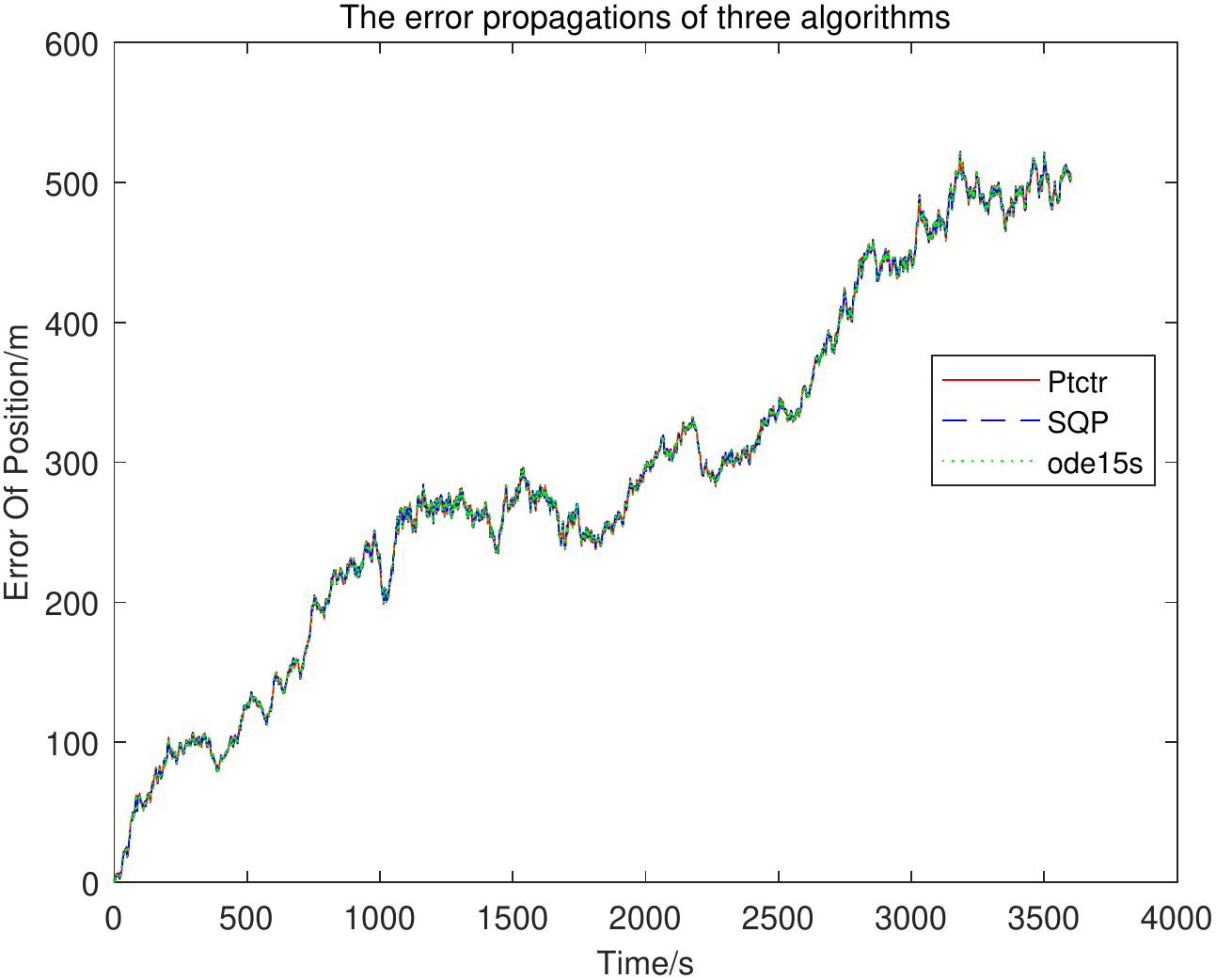}
      }
      \subfigure[trajectory 3]{
        \includegraphics[width=1\textwidth,height=0.25\textheight]{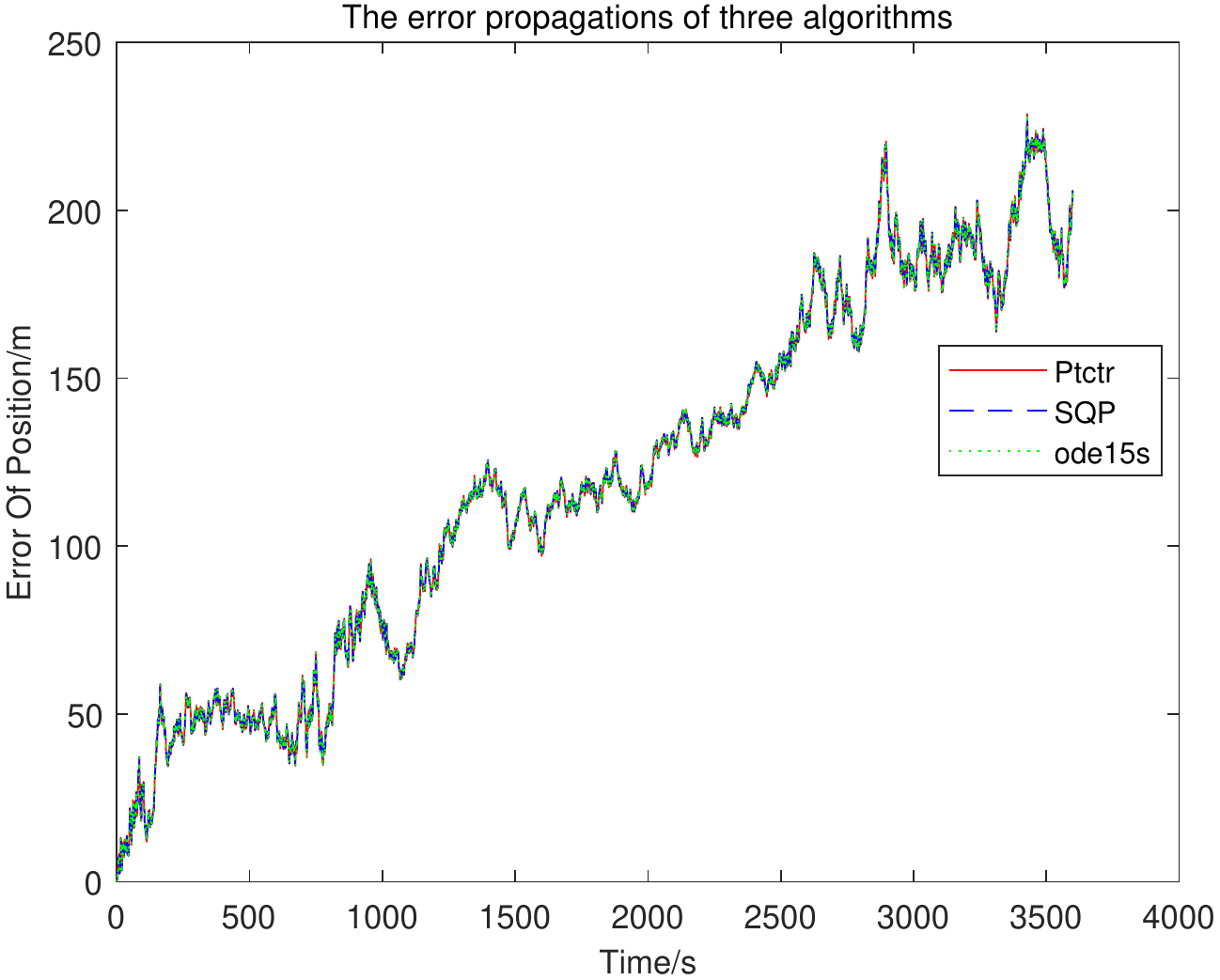}
      }
      \caption{The visual-inertial positioning problems
      with random errors.} \label{FIGWITHERR}
\end{minipage}
\end{figure}

First, we simulate three trajectories defined by equations
\eqref{TRA1}-\eqref{TRA3} without measurement errors via using Ptctr, SQP and
ode15s respectively to solve a series of optimization problems
\eqref{ERRDISFUN}-\eqref{ERRCON}. The initial point of the $(k+1)$-th
optimization problem is set by the optimal solution of the $k$-th optimization
problem. The simulation results are presented in Table \ref{TABCPUVIP} and
Figure \ref{FIGNOERR}. From Table \ref{TABCPUVIP} and Figure \ref{FIGNOERR},
we find the consumed time of Ptctr is about one-fifth of that ot SQP and
ode15s, and the trajectories computed by Ptctr and ode15s are more accurate
than those computed by SQP.

\vskip 2mm

Generally, the measurement data contain errors. According to the provided
industrial parameters, we assume that the measurement error $\varepsilon_{h}$
of altitude satisfies the Gaussian distribution $N(\mu_{h}, \; \sigma_{h})$ with the
variance $\sigma^{2}_{h} = 1$ and the mean $\mu_{h} = 0$. We assume that the
measurement error $\varepsilon_{d}$ of inertial distance satisfies an uniform
distribution on the interval $[-2.57, \; 2.57]$, i.e. $dist_{hor,\varepsilon}
= dist_{hor} + \varepsilon_{d}$.

\vskip 2mm

The angular error $\varepsilon_{\theta}$ of imaging affects the image coordinate
$\left(x_{pn}^{k}(\varepsilon_{\theta}), \; y_{pn}^{k}(\varepsilon_{\theta})\right)$
of the $n$-th landmark $(x_{ln}, \, y_{ln}, \, z_{ln})_{k}$ in the $k$-th camera
coordinate system as follows:
\begin{align}
      & \theta_{xn}^{k} =  \arctan \left(\frac{x_{k} - x_{ln}}{h_{n}^{k}}\right), \;
      \theta_{yn}^{k} = \arctan \left(\frac{y_{k} - y_{ln}}{h_{n}^{k}} \right), \nonumber \\
      & \theta_{xn}^{k}(\varepsilon_{\theta}) = \theta_{xn}^{k} + \varepsilon_{\theta}, \;
      \theta_{yn}^{k}(\varepsilon_{\theta}) = \theta_{yn}^{k} + \varepsilon_{\theta}, \nonumber \\
      & x_{pn}^{k}(\varepsilon_{\theta}) = f_{c} \tan \left(\theta_{xn}^{k}(\varepsilon_{\theta})\right), \;
      y_{pn}^{k}(\varepsilon_{\theta}) = f_{c} \tan \left(\theta_{yn}^{k}(\varepsilon_{\theta})\right),
    \; n = 1, \, 2, \ldots, N,    \label{PIXLANDMKERR}
\end{align}
where the altitude difference $h_{n}^{k}$ between the $k$-th camera and the $n$-th landmark
in the world coordinate system is computed by $h_{n}^{k} = z_{k} - z_{ln}$.
Here, we assume that the angular error $\varepsilon_{\theta}$ satisfies an uniform
distribution on the interval $[-0.2, \; 0.2]$.

\vskip 2mm

For conforming to the real environment of the unmanned aerial vehicle maintaining
the horizontal flight, we simulate three trajectories defined by equations
\eqref{TRA1}-\eqref{TRA3} with measurement errors via using Ptctr, SQP and
ode15s respectively to solve a series of optimization problems
\eqref{ERRDISFUN}-\eqref{ERRCON}. The measurement errors $\varepsilon_{h}, \;
\varepsilon_{d}, \; \varepsilon_{\theta}$ are stated by the previous several
paragraphs of this subsection. The numerical results are presented in Table \ref{TABCPUVIPERR}
and Figure \ref{FIGWITHERR}. From Table \ref{TABCPUVIPERR} and Figure \ref{FIGWITHERR},
we find that these three algorithms all work well for this problem and the consumed
time of Ptctr is one fifth of that of SQP and ode15s, respectively.

\vskip 2mm

\begin{table}[http]
  \caption{The visual-inertial positioning \\ problems without random errors.}
  \centering
  \begin{tabular}{llll}
    \hline
    \multirow{2}{*}{Trajectories} & Ptctr & SQP & ode15s  \\ \cline{2-4}
          & CPU Time (s)   & CPU Time (s) & CPU Time (s) \\ \hline
    trj 1 & 6.908 & 36.126 & 36.134  \\
    trj 2 & 6.891 & 37.608 & 37.614  \\
    trj 3 & 6.792 & 38.814  & 38.821\\ \hline
  \end{tabular} \label{TABCPUVIP}
\end{table}

\vskip 2mm

\begin{table}[http]
  \caption{The visual-inertial positioning \\ problems with random errors.}
  \centering
  \begin{tabular}{llll}
    \hline
    \multirow{2}{*}{Trajectories} & Ptctr & SQP & ode15s  \\ \cline{2-4}
          & CPU Time (s)   & CPU Time (s) & CPU Time (s) \\ \hline
    trj 1 & 6.572 & 37.885 & 37.893  \\
    trj 2 & 6.695 & 38.216 & 38.225  \\
    trj 3 & 5.618 & 32.023 & 32.029 \\ \hline
  \end{tabular} \label{TABCPUVIPERR}
\end{table}

\vskip 2mm

\section{Conclusion and Future Work}

\vskip 2mm

In this paper, we give a continuation method with the trusty time-stepping
scheme (Ptctr) for linearly equality-constrained optimization problems. Ptctr
only needs to solve a linear system of equations at every iteration, other
than the traditional optimization method such as SQP, which needs to solve
a quadratic programming subproblem at every iteration. This means that Ptctr
can save much more computational time than SQP. Numerical results also show that the
consumed time of Ptctr is about one fifth of that of SQP, which is the best
method of the three methods (SQP, PFM, ode15s). Furthermore, Ptctr
works well for simulating the trajectory of the unmanned aerial vehicle
maintaining the horizontal flight for a long time. Therefore, Ptctr is
worth investigating further, and we will extend it to the general nonlinear
optimization problem in the future.

\section*{Acknowledgments} This work was supported in part by Grant 61876199 from National
  Natural Science Foundation of China, Grant YBWL2011085 from Huawei Technologies
  Co., Ltd., and Grant YJCB2011003HI from the Innovation Research Program of Huawei
  Technologies Co., Ltd.. The authors are grateful to Prof. Hongchao Zhang, Prof. Li-Zhi Liao
  and two anonymous referees for their comments and suggestions which greatly
  improve the presentation of this paper.

\begin{appendix}

\section{Test Problems}

\noindent \textbf{Example 1.}
\begin{align}
   \quad  m & = n/2 \nonumber \\
   \min_{x \in \Re^{n}} \;  f(x) & = \sum_{k=1}^{n/2} \;\left(x_{2k-1}^{2} + 10x_{2k}^{2}\right), \;
   \text{subject to} \;  x_{2i-1} + x_{2i} = 4, \; i = 1, \, 2, \ldots, \, m. \nonumber
\end{align}
This problem is extended from the problem of \cite{Kim2010}. We assume that the
feasible initial point is $(2, \, 2, \, \ldots, \, 2, \, 2)$.

\vskip 2mm

\noindent \textbf{Example 2.}
\begin{align}
    \quad   m &= n/3 \nonumber \\
    \min_{x \in \Re^{n}} \;  f(x) & = \sum_{k=1}^{n/2} \; \left(\left(x_{2k-1} -2\right)^{2}
    + 2\left(x_{2k} - 1 \right)^{4}\right) - 5, \;
    \text{subject to} \; x_{3i-2} + 4x_{3i-1}+2x_{3i} = 3, \;  i = 1, \, 2, \ldots, \, n/3.
    \nonumber
\end{align}
We assume that the infeasible initial point is $(-0.5, \, 1.5, \, 1, \, 0, \, \ldots, \, 0, \, 0)$.

\vskip 2mm

\noindent \textbf{Example 3.}
\begin{align}
   \quad  m & = (2/3)n \nonumber \\
   \min_{x \in \Re^{n}} \;  f(x) & = \sum_{k=1}^{n}\; x_{k}^{2}, \;
   \text{subject to} \; x_{3i-2} + 2x_{3i-1} + x_{3i} = 1, \;
     2 x_{3i-2} - x_{3i-1} - 3 x_{3i} = 4, \; i = 1, \, 2, \ldots, \, n/3.  \nonumber
\end{align}
This problem is extended from the problem of \cite{Osborne2016}. The infeasible
initial point is $(1, \, 0.5, \, -1, \, \ldots, \, 1, \, 0.5, \, -1)$.

\vskip 2mm

\noindent \textbf{Example 4.}
\begin{align}
   \quad   m &= n/2 \nonumber \\
   \min_{x \in \Re^{n}} \; f(x) & = \sum_{k=1}^{n/2}\;\left(x_{2k-1}^{2} + x_{2k}^{6}\right) - 1,  \;
   \text{subject to} \;  x_{2i-1} + x_{2i} = 1, \; i = 1, \, 2, \, \ldots, \, n/2. \nonumber
\end{align}
This problem is modified from the problem of \cite{MAK2019}. We assume that the
infeasible initial point is $(1, \, 1, \, \ldots, \, 1)$.

\vskip 2mm

\noindent \textbf{Example 5.}
\begin{align}
   \quad  m & = n/2 \nonumber \\
   \min_{x \in \Re^{n}} \;  f(x) & = \sum_{k=1}^{n/2}\;\left(\left(x_{2k-1} -2\right)^{4}
   + 2\left(x_{2k} -1\right)^{6}\right) - 5, \;
    \text{subject to} \; x_{2i-1} + 4x_{2i} = 3, \; i = 1, \, 2, \, \ldots, \, m. \nonumber
\end{align}
We assume that the feasible initial point is $(-1, \, 1,\, -1, \, 1, \, \ldots, \, -1, \, 1)$.

\vskip 2mm

\noindent \textbf{Example 6.}
\begin{align}
   \quad m & = (2/3)n \nonumber \\
   \min_{x \in \Re^{n}} \;  f(x) &=
   \sum_{k=1}^{n/3}\;\left(x_{3k-2}^{2} + x_{3k-1}^{4} + x_{3k}^{6}\right), \nonumber \\
    \text{subject to} \; &  {x_{3i-2}} + 2{x_{3i-1}} + {x_{3i}} = 1, \;
      2{x_{3i-2}} - {x_{3i-1}} - 3{x_{3i}} = 4, \; i = 1, \, 2, \, \ldots, \, m/2. \nonumber
\end{align}
This problem is extended from the problem of \cite{Osborne2016}. We assume that the
infeasible initial point is $(2, \, 0, \, \ldots, \, 0)$.

\vskip 2mm

\noindent \textbf{Example 7.}
\begin{align}
   \quad m &= n/2 \nonumber \\
   \min_{x \in \Re^{n}} \;  f(x) & = \sum_{k=1}^{n/2}\;\left(x_{2k-1}^{4} + 3x_{2k}^{2}\right), \;
   \text{subject to} \;  x_{2i-1} + x_{2i} = 4, \; i = 1, \, 2, \, \ldots, \, n/2.  \nonumber
\end{align}
This problem is extended from the problem of \cite{Carlberg2009}. We assume that the infeasible initial
point is $(2, \, 2, \, 0, \, \ldots, \, 0, \, 0)$.

\vskip 2mm

\noindent \textbf{Example 8.}
\begin{align}
   \quad  m & = n/3 \nonumber \\
   \min_{x \in \Re{^n}} \;  f(x) & = \sum_{k=1}^{n/3}\;\left(x_{3k-2}^{2} + x_{3k-2}^{2} \, x_{3k}^{2}
   + 2x_{3k-2} \, x_{3k-1} + x_{3k-1}^{4} + 8x_{3k-1}\right), \nonumber \\
   \text{subject to} \; & 2x_{3i-2} + 5x_{3i-1}+x_{3i} = 3, \; i = 1, \, 2, \, \ldots, \, m. \nonumber
\end{align}
We assume that the infeasible initial point is $(1.5, \, 0, \, 0, \, \ldots, \, 0)$.

\vskip 2mm

\noindent \textbf{Example 9.}
\begin{align}
   \quad m &= n/2 \nonumber \\
   \min_{x \in \Re^{n}} \;  f(x) & = \sum_{k =1}^{n/2} \; \left(x_{2k-1}^{4} + 10x_{2k}^{6}\right),
   \;
   \text{subject to} \; x_{2i-1} + x_{2i} =4, \; i = 1, \, 2, \, \ldots, \, m.  \nonumber
\end{align}
This problem is extended from the problem of \cite{Kim2010}. We assume that the feasible
initial point is $(2, \, 2, \, \ldots, \, 2, \, 2)$.

\vskip 2mm

\noindent \textbf{Example 10.}
\begin{align}
   \quad m & = n/3 \nonumber \\
   \min_{x \in \Re^{n}} \;  f(x) & = \sum_{k=1}^{n/3}\;\left(x_{3k-2}^{8} + x_{3k-1}^{6} + x_{3k}^{2}\, \right),
   \;
   \text{subject to} \;  x_{3i-2} + 2x_{3i-1} + 2x_{3i} =1, \; i = 1, \, 2, \, \ldots, \,m. \nonumber
\end{align}
This problem is modified from the problem of \cite{Yamashita1980}.
The feasible initial point is $(1, \, 0, \, 0, \, \ldots, \, 1, \, 0, \, 0)$.

\vskip 2mm

\end{appendix}

\end{document}